\documentclass[american,english]{amsart}
\usepackage[latin1]{inputenc}
\usepackage{babel}
\usepackage{graphics}
\usepackage{amssymb}

\makeatletter

\providecommand{\LyX}{L\kern-.1667em\lower.25em\hbox{Y}\kern-.125emX\@}
\newcommand{\noun}[1]{\textsc{#1}}

 \theoremstyle{plain}    
 \newtheorem{thm}{Theorem}[section]
 \numberwithin{equation}{section} 
 \numberwithin{figure}{section} 
 \theoremstyle{plain}    
 \newtheorem{cor}[thm]{Corollary} 
 \theoremstyle{plain}    
 \newtheorem{lem}[thm]{Lemma} 

 \theoremstyle{plain}    
 \newtheorem{obs}[thm]{Observation}

 \theoremstyle{plain}

 \theoremstyle{definition}
 \newtheorem{defn}[thm]{Definition}
 \theoremstyle{definition}
 \newtheorem{asmp}[thm]{Assumption}
 \theoremstyle{definition}
 \newtheorem{conv}[thm]{Convention}
 \theoremstyle{remark}
 \newtheorem{rem}[thm]{Remark}
 \theoremstyle{remark}    
 \newtheorem{notation}[thm]{Notation} 

\makeatother
\begin{document}

\title{a dichotomy in classifying quantifiers for finite models }

\author{saharon shelah${}^{1}$ }

\thanks{${}^{1}$This research was supported by The Israel Science Foundation.
   Publication 801.}

\author{mor doron}

\begin{abstract}
We consider a family \( \frak {U} \) of finite universes. The second
order quantifier \( Q_{\frak {R}} \), means for each \( U\in \frak {U} \)
quantifying over a set of \( n(\frak {R}) \)-place relations isomorphic
esto a given relation. We define a natural partial order on such quantifiers
called interpretability. We show that for every \( Q_{\frak {R}} \),
ever \( Q_{\frak {R}} \) is interpretable by quantifying over subsets
of \( U \) and one to one functions on \( U \) both of bounded order,
or the logic \( L(Q_{\frak {R}}) \) (first order logic plus the quantifier
\( Q_{\frak {R}} \)) is undecidable. 
\end{abstract}
\maketitle
\tableofcontents{}

\section{Introduction}

\subsection{Background}

In this work we continue \cite{key-1}, but it is self contained and
the reader may read it independently. Our aim is to analyze and classify
second order quantifiers in finite model theory. The quantifiers will
be defined as follows:\begin{itemize}
\item[(*)]

Let \( U \) be a finite universe, and \( n \) a natural number.
Let \( K \) be a class of \( n \)-place relations on \( U \) closed
under permutations of \( U \). Define \( Q_{K} \) to be the quantifier
ranging over the relations in \( K \).

\end{itemize}

We will usually work on quantifiers of the form \( Q_{R}=Q_{K_{R}} \)
where \( R \) is a \( n \)-place relation over \( U \) and \( K_{R} \)
is defined by: \( K_{R}:=\{R'\subseteq {}^{n}U:(U,R)\approx (U,R')\} \).
We define below two partial orders on the class of such quantifiers,
called: interpretability and expressibility. It will be interesting
to consider the class \( K \) of \( n \)-place relations definable
in some logic \( \frak {L} \), that is such that there exists a formula
\( \varphi (r)\in \frak {L} \) (\( r \) is a \( n \)-place relation
symbol) and \( R\in K \) iff \( (U,R)\models \varphi (r) \). In
\cite{key-2} the problem was solved for the case: \( K \) is definable
in first order logic and \( U \) is infinite. It was shown that in
this case \( Q_{K} \) is equivalent (in the sense of interpretability)
to one of only four quantifiers: trivial (first order), monadic, quantifying
over 1-1 functions or full second order. A revue paper of this result
is \cite{key-6}. If we do not assume \( Q_{K} \) to be first order
definable but keep assuming \( U \) is infinite we get a classification
of \( Q_{K} \) by equivalence relations. Formally from \cite{key-3}
we have: 

\begin{thm}
\label{Infinite} Let \( U \) be an infinite countable universe,
and \( K \) be as in \( (*) \). Then there exist a family \( \bf {E} \)
of equivalence relations on \( U \), such that \( Q_{K} \) and \( Q_{\bf {E}} \)
are equivalent (each is interpretable by the other). 
\end{thm}
We remark that if \( U \) is infinite not nessesarily countable then
the situation is more complicated, but if we assume \( L=V \) then
we have the same result. \cite{key-1} deals with the case \( U \)
is finite. Under this assumption we get a reasonable understanding
of \( Q_{R} \), we can {}``bound{}'' it between two simple and
close quantifiers (close meaning that the size of one is a polynomial
in the size of the other). Formally: 

\begin{thm}
\label{Finite}Let \( U \) be a finite universe, and \( R \) a \( n \)-place
relation on \( U \). Then there exist a natural number \( \lambda =\lambda (R) \),
and equivalence relation \( E \) on \( U \) such that uniformly
we have: 
\begin{enumerate}
\item \( Q_{E} \) and \( Q^{1-1}_{\lambda } \) are interpretable by \( Q_{R} \)
(\( Q^{1-1}_{\lambda } \) is the quantifier ranging over 1-1 partial
functions with domain \( \leq \lambda  \). 
\item If \( |U|\geq \lambda ^{n} \) then \( Q_{R} \) is expressible by
\( \{Q_{E},Q^{1-1}_{\lambda ^{n}}\} \).
\item If \( |U|<\lambda ^{n} \) then every \( 2 \)-place relation on a
subset \( A\subseteq U \) with cardinality \( \leq |U|^{1/2n} \)
is interpretable by \( Q_{R} \).
\end{enumerate}
Where {}``uniformly{}'' means the formulas used to express and interpret
are independent of \( U \) and depend on \( n \) alone. 
\end{thm}
In case (2) of the theorem if we want to have {}``interpretable{}''
instead of {}``expressible{}'' then the situation is more complicated
and we deal with it in this paper. Since \( U \) is a {}``large{}''
universe we check the {}``asymptotic behavior{}'', that is we consider
a class \( \frak {U} \) of finite universes with unbounded cardinality.
For each \( U\in \frak {U} \) let \( \frak {R}[U]\subseteq {}^{n}U \)
be an \( n \)-place relation on \( U \). We will see that there
is a dichotomy y in the behavior of \( Q_{\frak {R}[U]} \), that relates
to cases (1) and (2) of theorem \ref{Finite}. Formally we prove: 

\begin{thm}
\label{theorem}Let \( \frak {R} \) be as above. Then exactly one
of the following conditions holds:
\begin{enumerate}
\item \( Q_{\frak {R}[U]} \) is uniformly interpretable by 1-1 functions
and \( 1 \)-place relations both of bounded cardinality.
\item For each \( m\in \mathbb {N} \), there exist \( U\in \frak {U} \)
such that we can uniformly interpret number theory up to \( m \),
by \( Q_{\frak {R}[U]} \).
\end{enumerate}
\end{thm}
We prove this theorem in sections 3 to 6. In section 3 we analyze
the situation, and give a condition for the dichotomy. In section
4 we prove that if the condition of section 3 hold then part (2) of
theorem \ref{theorem} is satisfied. In section 5 we prove, for the
\( 2 \)-place case that if the condition does not hold then part
(1) of the theorem is satisfied. In section 6 we prove the same for
the \( n \)-place case. In section 2 we show that in the finite case
we can not get a full understanding of \( Q_{\frak {K}} \) similar
to what we have in the countable case (not even for expressibility).

\subsection{Notations Conventions And Primary Definitions}

\begin{conv}
\label{Motation}\(  \)
\begin{enumerate}
\item \( \frak {U} \) is a class of finite universes, possibly with repetitions.
So formally: \( \frak {U}=\{U_{i}:i\in \frak {I}\} \) for an index
class \( \frak {I} \) and we allow \( U_{i}=U_{j} \) for \( i\neq j\in \frak {I} \).
We will usually not be so formal and will write \( U\in \frak {U} \)
and it should be understood as \( i\in \frak {I} \) and \( U=U_{i} \).
We assume \( sup\{|U|:U\in \frak {U}\}=\aleph _{0} \).
\item \( \frak {K} \) is a function on \( \frak {U} \) and for all \( U\in \frak {U} \),
\( \frak {K}[U] \) is a set of \( n \)-place relations on \( U \)
(where \( n=n(\frak {K}) \) is a natural number), closed under permutations
of \( U \). This means: if \( R_{1},R_{2}\subseteq {}^{n}U \) and
\( (U,R_{1})\approx (U,R_{2}) \) then \( R_{1}\in \frak {K}[U]\Leftrightarrow R_{2}\in \frak {K}[U] \).
\item \( \overline{\frak {K}} \) is a sequence of such functions. We write
\( \overline{\frak {K}}=(\frak {K}_{0},...,\frak {K}_{lg(\overline{\frak {K}})-1}) \).
\item \( \frak {R} \) is a function on \( \frak {U} \) and for each \( U\in \frak {U} \),
\( \frak {R}[U] \) is a \( n \)-place relation over \( U \) (where
\( n=n(\frak {R}) \) is a natural number).
\item \( r \) is a \( n(\frak {R}) \)-place relation symbol.
\item For all \( U\in \frak {U} \) if \( S \) is a \( n \)-place relation
on \( U \), and \( F \) is a \( m \)-place function on \( U \),
then \( s \) and \( f \) are a \( n \)-place relation symbol and
a \( m \)-place function symbol respectively. We write \( (U,S)\models s(\overline{a}) \)
iff \( \overline{a}\in s \), and \( (U,F)\models f(\overline{b})=c \)
iff \( F(\overline{b})=c \). (That is for all \( c\in U,\bar{a}\in {}^{n}U,\bar{b}\in {}^{m}U \)). 
\item For all \( U\in \frak {U} \) and \( n\in \omega  \), \( \bar{a}\in {}^{n}U \)
is a sequence of \( n \) elements in \( U \). We write: \( \bar{a}=(a_{0},...,a_{n-1}) \),
and \( lg(\bar{a})=n \).
\end{enumerate}
\end{conv}
\begin{defn}
For all \( \frak {K} \) as in \ref{Motation}.2 we define the second
order quantifier \( Q_{\frak {K}} \) to range over all relations
in \( \frak {K} \). Formally we define the logic \( L(Q_{\frak {K}_{1}},...,Q_{\frak {K}_{m}}) \)
to be first order logic but we allow formulas of the form \( (Q_{\frak {K}_{i}}r)\varphi (r) \)
(\( r \) is a \( n(\frak {K}_{i}) \)-place relation symbol) for
all \( 1\leq i\leq m \). Satisfaction is defined only for models
with universe \( U\in \frak {U} \) as follows: \( \models (Q_{\frak {K}_{i}}r)\varphi (r) \)
iff there exists \( R^{0}\in \frak {K}_{i}[U] \) such that \( (U,R^{0})\models \varphi (r) \).
\begin{defn}
\label{Definable}We say that \( \frak {K} \) (or \( Q_{\frak {K}} \))
is definable in some logic \( \frak {L} \) iff there exists a formula
\( \varphi (r)\in \frak {L} \) (\( r \) is a \( n(\frak {K}) \)-place
relation symbol) such that for all \( U\in \frak {U} \) and \( R\subseteq {}^{n(\frak {K})}U \):\[
(U,R)\models \varphi (r)\Longleftrightarrow R\in \frak {K}[U]\]

\end{defn}
\end{defn}
\begin{notation}
For \( \frak {R} \) as in \ref{Motation}.4 we note \( Q_{\frak {R}} \)
by \( Q_{\frak {K}_{\frak {R}}} \) where \( \frak {K}=\frak {K}_{\frak {R}} \)
is defined by:\[
\frak {K}[U]:=\{R^{1}\subseteq {}^{n(\frak {R})}U:(U,R^{1})\approx (U,\frak {R}[U])\}\]

\end{notation}
\begin{defn}
\label{int exp}\(  \)
\begin{enumerate}
\item We say that \( Q_{\frak {K}_{1}} \) is interpretable by \( Q_{\frak {K}_{2}} \)
and write \( Q_{\frak {K}_{1}}\leq _{int}Q_{\frak {K}_{2}} \) if
there exist \( k^{*}\in \omega  \) and first order formulas: \( \varphi _{k}(\overline{x},\overline{r})=\varphi _{k}(x_{0},...,x_{n(\frak {K}_{1})-1},r_{0},...,r_{m-1}) \)
for \( k<k^{*} \) (each \( r_{l} \) is a \( n(\frak {K}_{2}) \)-place
relation symbol) and the following holds: \\
\begin{itemize}
\item[(*)] For all \( U\in \frak {U} \) and \( R\in \frak {K}_{1}[U] \)
there exists \( k<k^{*} \) and \( R_{0},...,R_{m-1}\in \frak {K}_{2}[U] \)
such that \( (U,R_{0},...,R_{m-1})\models (\forall \overline{x})[R(\overline{x})\equiv \varphi _{k}(\overline{x},r_{0},...,r_{m-1})] \).\\
\end {itemize}
\item We say that \( Q_{\frak {K}_{1}} \) is expressible by \( Q_{\frak {K}_{2}} \)
and write \( Q_{\frak {K}_{1}}\leq _{exp}Q_{\frak {K}_{2}} \) if
there exist \( k^{*}\in \omega  \) and formulas in the logic \( L(Q_{\frak {R}_{2}}) \):
\( \varphi _{k}(\overline{x},\overline{r})=\varphi _{k}(x_{0},...,x_{n(\frak {K}_{1})-1},r_{0},...,r_{m-1}) \)
for \( k<k^{*} \) (each \( r_{l} \) is a \( n(\frak {K}_{2}) \)-place
relation symbol) and \( (*) \) holds. 
\item In (1) and (2) if \( k^{*}=1 \) we write \( Q_{\frak {K}_{1}}\leq _{1-int}Q_{\frak {K}_{2}} \)
and \( Q_{\frak {K}_{1}}\leq _{1-exp}Q_{\frak {K}_{2}} \) respectively.
\item We write \( Q_{\frak {K}_{1}}\equiv _{int}Q_{\frak {K}_{2}} \) if
\( Q_{\frak {K}_{1}}\leq _{int}Q_{\frak {K}_{2}} \) and \( Q_{\frak {K}_{2}}\leq _{int}Q_{\frak {K}_{1}} \).
\( \equiv _{exp} \) is defined in the same way.
\item We define \( Q_{\frak {K}}\leq _{int}\{Q_{\frak {K}_{0}},...,Q_{\frak {K}_{l-1}}\} \)
as in (1) only in \( (*) \) we allow \( R_{0},...,R_{m-1}\in \bigcup _{i=0}^{l-1}\frak {K}_{i}[U] \).
We write \( Q_{\overline{\frak {K}}}=\{Q_{\frak {K}_{0}},...,Q_{\frak {K}_{lg(\overline{\frak {K}})-1}}\} \)
when \( \overline{\frak {K}}=(\frak {K}_{0},...,\frak {K}_{lg(\overline{\frak {K}})-1}) \).
In the same way we define for \( \leq _{exp} \).
\item We define \( Q_{\overline{\frak {K}^{1}}}\leq _{int}Q_{\overline{\frak {K}^{2}}} \)
if \( Q_{\frak {K}^{1}_{i}}\leq _{int}Q_{\overline{\frak {K}^{2}}} \)
for all \( i<lg(\overline{\frak {K}^{1}}) \) again when \( \overline{\frak {K}^{1}}=(\frak {K}^{1}_{0},...,\frak {K}^{1}_{lg(\overline{\frak {K}^{1}})-1}) \).
In the same way we define for \( \leq _{exp} \).
\end{enumerate}
\end{defn}
\begin{lem}
\label{hirarchy}\(  \) 
\begin{enumerate}
\item \( \leq _{int} \) and \( \leq _{exp} \) are partial orders, Hence
\( \equiv _{int} \) and \( \equiv _{exp} \) are equivalence relations
on the class of quantifiers of the form \( Q_{\frak {K}} \).
\item \( Q_{\overline{\frak {K}^{1}}}\leq _{int}Q_{\overline{\frak {K}^{2}}} \)
implies \( Q_{\overline{\frak {K}^{1}}}\leq _{exp}Q_{\overline{\frak {K}^{2}}} \).
\end{enumerate}
\end{lem}
\begin{proof}
Straight.
\end{proof}
So \( \leq _{exp} \) gives a hierarchy on on logics of the form \( \frak {L}(Q_{\overline{\frak {K}}}) \),
i.e under the assumptions of lemma \ref{hirarchy} the expressive
power of \( \frak {L}(Q_{\overline{\frak {K}^{2}}}) \) is at least
as strong as that of \( \frak {L}(Q_{\overline{\frak {K}^{1}}}) \).

\begin{lem}
Let \( \frak {L} \) be some logic and assume \( \overline{\frak {K}^{1}},\overline{\frak {K}^{2}} \)
are definable in \( \frak {L} \) (that is every \( \frak {K}_{i}^{l} \)
is, see definition \ref{Definable}) and \( Q_{\overline{\frak {K}^{1}}}\leq _{exp}Q_{\overline{\frak {K}^{2}}} \)
then:
\begin{enumerate}
\item there exists a computable function that attach to every formula in
\( \frak {L}(Q_{\overline{\frak {K}^{1}}}) \) an equivalent formula
in \( \frak {L}(Q_{\overline{\frak {K}^{2}}}) \).
\item the set of valid sentences in \( \frak {L}(Q_{\overline{\frak {K}^{1}}}) \)
is recursive from the set of valid sentences in \( \frak {L}(Q_{\overline{\frak {K}^{2}}}) \).
\end{enumerate}
\end{lem}
\begin{proof}
Easy.
\end{proof}

\subsection{Summation of Previous Results.}

We will use the following results. Proofs can be found in \cite{key-1}.

\begin{defn}
\label{Q}\(  \)
\begin{enumerate}
\item let \( \lambda  \) be a function from \( \frak {U} \) to \( \mathbb {N} \)
such that \( \lambda [U]\leq |U|/2 \). Define \( \frak {K}^{mon}_{\lambda } \)
by \( \frak {K}_{\lambda }^{mon}[U]:=\{A\subseteq U:|A|=\lambda [U]\} \).
We note \( Q_{\frak {K}^{mon}_{\lambda }} \) by \( Q_{\lambda }^{mon} \).
\item For \( \lambda  \) as above define \( \frak {K}^{mon}_{\leq \lambda } \)
by \( \frak {K}^{mon}_{\leq \lambda }[U]:=\bigcup \{\frak {K}^{mon}_{\mu }:\mu \leq \lambda \} \).
We note \( Q_{\frak {K}^{mon}_{\leq \lambda }} \) by \( Q^{mon}_{\leq \lambda } \).
\item For \( \lambda  \) as above define \( \frak {K}^{1-1}_{\lambda } \)
by \( \frak {K}_{\lambda }^{1-1}[U]:=\{f:U\rightarrow U:|Dom(f)|=\lambda [U],\, f\, one\, to\, one\} \).
We note \( Q_{\frak {K}^{1-1}_{\lambda }} \) by \( Q_{\lambda }^{1-1} \).
\item For \( \lambda  \) as above define \( \frak {K}^{1-1}_{\leq \lambda } \)
by \( \frak {K}^{1-1}_{\leq \lambda }[U]:=\bigcup \{\frak {K}^{1-1}_{\mu }:\mu \leq \lambda \} \).
We note \( Q_{\frak {K}^{1-1}_{\leq \lambda }} \) by \( Q^{1-1}_{\leq \lambda } \).
\item Let \( \lambda  \) and \( \mu  \) be functions from \( \frak {U} \)
to \( \mathbb {N} \). Define \( \frak {K}^{eq}_{\lambda ,\mu } \)
as follows: \( \frak {K}^{eq}_{\lambda ,\mu }[U] \) is the collection
of all equivalence relations on subsets of \( U \) with exactly \( \lambda [U] \)
classes, and the size of each class is \( \mu [U] \). We note \( Q_{\frak {K}^{eq}_{\lambda ,\mu }} \)
by \( Q_{\lambda ,\mu }^{eq} \). 
\item Let \( \lambda  \) and \( \mu  \) be as in 5. Define \( \frak {K}^{eq}_{\leq \lambda ,\leq \mu } \)
as follows: \( \frak {K}^{eq}_{\leq \lambda ,\leq \mu } \) is the
collection of all equivalence relations on subsets of \( U \) with
at most \( \lambda [U] \) classes, and the size of each at most is
\( \mu [U] \). We note \( Q_{\frak {K}^{eq}_{\leq \lambda ,\leq \mu }} \)
by \( Q^{eq}_{\leq \lambda ,\leq \mu } \). 
\end{enumerate}
\end{defn}
\begin{rem}
In \cite{key-1} it is proved that \( Q^{mon}_{\lambda }\equiv _{int}Q^{mon}_{\leq \lambda } \)
and \( Q^{1-1}_{\lambda }\equiv _{int}Q^{1-1}_{\leq \lambda } \)
so we will usually not distinguish between them.
\end{rem}
\begin{lem}
\label{Int k,k}Let \( \lambda  \) be a function from \( \frak {U} \)
to \( \mathbb {N} \), and \( \frak {E} \) a 2-place relation on
\( \frak {U} \) such that for all \( U\in \frak {U} \), \( \frak {E}[U] \)
is an equivalence relation with at least \( \lambda [U] \) classes
each of which has at least \( \lambda [U] \) elements (and possibly
smaller classes). Then \( Q_{\lambda ,\lambda }\leq _{1-int}Q_{\frak {E}} \).
\end{lem}
\begin{proof}
straight foreword. The interpreting formula is \( \varphi (x,y,s_{0},s_{1},s_{2}):=s_{0}(x,y)\wedge \lnot s_{1}(x,y)\wedge s_{2}(x,y) \).
(See \cite{key-1} for similar proofs).
\end{proof}
\begin{thm}
\label{Lambda_0}For every \( \frak {R} \) as in \ref{Motation}.4
there exists a function \( \lambda _{0}=\lambda _{0}(\frak {R}) \)
from \( \frak {U} \) to \( \mathbb {N} \) such that: 
\begin{enumerate}
\item \( Q^{mon}_{\lambda _{0}}\leq _{int}Q_{\frak {R}} \). 
\item There exists \( \frak {R}_{1} \) with \( n=n(\frak {R})=n(\frak {R}_{1}) \)
and \( |Dom(\frak {R}_{1}[U])|\leq \lambda _{0}[U]+n \) for all \( U\in \frak {U} \),
such that \( Q_{\frak {R}}\equiv _{int}\{Q_{\frak {R}_{1}},Q_{\lambda _{0}}^{mon}\} \). 
\end{enumerate}
The interpretation is done uniformly, that is the formulas used are
independent of \( \frak {R} \) (depend on \( n(\frak {R}) \) alone).
\begin{thm}
\label{Lambda_1}For each \( \frak {R} \) there exists a function
\( \lambda _{1}=\lambda _{1}(\frak {R}) \) from \( \frak {U} \)
to \( \mathbb {N} \) such that uniformly: \( Q_{\frak {R}}\equiv _{int}\{Q_{\lambda _{0}}^{mon},Q_{\lambda _{1}}^{1-1},Q_{\frak {R}_{1}},Q_{\frak {E}}\} \),
where \( n=n(\frak {R})=n(\frak {R}_{1}) \) and for all \( U\in \frak {U} \),
\( |Dom(\frak {R}_{1}[U])|\leq n\cdot \lambda _{1}[U] \) and \( \frak {E}[U] \)
is an equivalence relation on \( U \).
\end{thm}
\end{thm}
\begin{rem}
\label{chenge lambda}In the proof of theorem \ref{theorem} we can assume
without loss of generality that for all \( U\in \frak {U} \), \( |\frak {R}[U]|\leq \lambda _{0}[U]+n(\frak {R}) \),
this is true since we can interpret \( \frak {R}_{1} \) instead of
\( \frak {R} \) (see \ref{Lambda_0}). Similarly using \ref{Lambda_1}
we can assume \( |\frak {R}[U]|\leq \lambda _{1}[U]\cdot n(\frak {R}) \).
Here we have an equivalence relation \( \frak {E} \) that can change
the bounds but the change will not be significant. Note also that
\( Q^{1-1}_{\lambda _{1}}\equiv _{int}Q^{1-1}_{n\cdot \lambda _{1}} \)
(for all \( n\in \omega  \)). So in the simple case of the dichotomy
(theorem \ref{simple}) we prove \( Q_{\frak {R}}\leq _{int}\{Q^{mon}_{\lambda _{0}},Q^{1-1}_{\lambda _{1}}\} \)
but in the proof we will not pay attention to the size of the sets
and functions we use.
\end{rem}

\section{Limitations on The Classification of \protect\( Q_{\frak {K}}\protect \)
in The Finite}

In this section we show that unlike the countable case in which we
had an understanding of \( Q_{\frak {K}} \) by equivalence relations,
in the finite case there are classes of relations we can not express. 

\begin{defn}
For all \( n\in \omega  \) define \( \frak {K}_{n} \) by: \( \frak {K}_{n}[U]:=\{R:R\subseteq {}^{n}U\} \)
for all \( U\in \frak {U} \).
\end{defn}
\begin{obs}
For all \( n\in \omega  \): \( Q_{\frak {K}_{n+1}}\nleq _{exp}Q_{\frak {K}_{n}} \).
\end{obs}
\begin{proof}
Suppose \( Q_{\frak {K}_{n+1}}\leq _{exp}Q_{\frak {K}_{n}} \), and
assume that the formulas used for expressing are \( \varphi _{k}(\bar{x},r_{0},...,r_{m_{k}-1}) \)
for \( k<k^{*} \). Note \( m=max\{m_{k}:k<k^{*}\}\cup \{k^{*}\} \)
and let \( U\in \frak {U} \). Then by these formulas we can express
at most \( m^{2}\cdot |\frak {K}_{n}[U]| \) different relation. Since
\( |\frak {K}_{n}[U]|=2^{|U|^{n}} \), if we choose \( U \) such
that \( |U|>\root n\of {log_{2}(m^{2})} \) we get \( |U|^{n}(|U|-1)>log_{2}(m^{2}) \),
hence \( 2^{|U|^{n}(|U|-1)}>m^{2} \), and hence \( 2^{|U|^{n+1}}>m^{2}\cdot 2^{|U|^{n}} \).
So the maximal number of different expressible relations is smaller
than \( |\frak {K}_{n+1}[U]| \), a contradiction. 
\end{proof}
We have that for \( n>2 \), \( Q_{\frak {K}_{n}} \) is not expressible
by equivalence relations, unlike the countable case (see \ref{Infinite}).
Moreover we have:

\begin{obs}
For all \( n\geq 2 \):
\begin{enumerate}
\item \( Q_{\frak {K}_{n}}\not \leq _{exp}Q^{1-1} \).
\item \( Q_{\frak {K}_{n}}\not \leq _{exp}Q^{eq} \).
\end{enumerate}
\end{obs}
\begin{proof}
We prove (1). again suppose \( Q_{\frak {K}_{n}}\leq _{exp}Q^{1-1} \),
and we use the notations of the previous proof. Note that \( \frak {K}^{1-1}[U]=|U|! \),
and for \( |U| \) large enough we have \( |U|!<2^{|U|\cdot log(|U|)\cdot c} \)
where \( c \) is some constant. Moreover for all \( n\geq 2 \) and
\( |U| \) large enough we have \( |U|\cdot log(|U|)\cdot c<|U|^{n} \).
So we get: \( m\cdot |K^{1-1}[U]|<2^{|U|^{n}} \) which means the
number of relations expressible is smaller than \( |\frak {K}_{n}[U]| \),
a contradiction.

The proof of (2) is similar using: \( \frak {K}^{eq}[U]\leq |U|^{|U|}\leq 2^{|U|\cdot log(|U|)\cdot c} \). 
\end{proof}
We get that in the finite case even for \( n(\frak {K})=2 \), we
can not express every \( Q_{\frak {K}} \) by 1-1 functions and equivalence
relations.

\section{Primary Analyses}

\begin{asmp}
\relax From here on, unless said otherwise, we assume that \( \frak {R} \)
(see \ref{Motation}.4) is fixed and \( \lambda _{i}=\lambda _{i}(\frak {R}) \)
for \( i\in \{0,1\} \) (see \ref{Lambda_0} and \ref{Lambda_1}).
\end{asmp}
In this section we start the analyses of \( Q_{\frak {R}} \). For
each universe \( U \) we define a natural number \( k \) which is
the maximal size, in some sense, of an equivalence relation on \( U \)
interpretable by \( \frak {R}[U] \). The size of \( k \) is an indicator
of the degree of {}``complexity{}'' of \( \frak {R} \). We will
show that there is a dichotomy, ether \( \frak {R} \) is very {}``complex{}''
or it is {}``simple{}''. this is made precise below. 

\begin{defn}
Let \( \tau =\{f_{0},...,f_{m_{1}},s_{0},...,s_{m_{2}},c_{0},...,c_{m_{3}}\} \)
be a vocabulary, that is \( f_{i} \) are \( n(f_{i}) \)-place function
symbols, \( s_{i} \) are \( n(s_{i}) \)-place relation symbols and
\( c_{i} \) are individual constants. Define: 
\begin{enumerate}
\item for all \( U\in \frak {U} \) a model for \( \tau  \) on \( U \)
is \( M=(U,f_{0}^{M},...,f^{M}_{m_{1}},s^{M}_{0},...,s^{M}_{m_{2}},c^{M}_{0},...,c^{M}_{m_{3}}) \),
where \( f^{M}_{i} \) are \( n(f_{i}) \)-place functions on \( U \),
\( s_{i}^{M} \) are \( n(s_{i}) \)-place relations on U and \( c_{i}^{M}\in U \).
\( U \) is called the universe of \( M \) and noted by \( |M| \).
\item a model for \( \tau  \) on \( \frak {U} \) noted by \( \frak {M} \)
is a function from \( \frak {U} \) such that for all \( U\in \frak {U} \),
\( \frak {M}[U] \) is a model for \( \tau  \) on \( U \). Note
that the function \( U\mapsto (U,\frak {R}[U]) \) is a model for
\( \{r\} \) on \( \frak {U} \), we will not be as formal and say
that \( \frak {R} \) is. 
\item Assume \( r\in \tau  \). We say that \( \frak {M} \) expand (or
is an expansion of) \( \frak {R} \) if for all \( U\in \frak {U} \),
\( r^{\frak {M}[U]}=\frak {R}[U] \). More generally:
\item Let \( \tau \subseteq \tau ' \) be dictionaries, and let \( \frak {M} \)
and \( \frak {M}' \) be models on \( \frak {U} \) for \( \tau  \)
and \( \tau ' \) respectively. We say that \( \frak {M}' \) expand
\( \frak {M} \) if \( \frak {M}'|\tau =\frak {M} \). That means
for all \( U\in \frak {U} \) and \( f_{i}\in \tau  \), \( f_{i}^{\frak {M}'[U]}=f_{i}^{\frak {M}[U]} \),
and similarly for relation symbols and constants.
\item We call \( \tau  \) \noun{simple} if:

\begin{enumerate}
\item \( \tau  \) is finite.
\item For all \( i\leq m_{1} \), \( n(f_{i})=1 \).
\item For all \( i\leq m_{2} \), \( n(s_{i})=1 \).
\end{enumerate}
\item We call \( M \) a \noun{simple} model for \( \tau  \) on \( U \)
if:

\begin{enumerate}
\item \( \tau  \) is simple.
\item \( M \) is a model for \( \tau  \) on \( U \).
\item For all \( i\leq m_{1} \), \( f_{i}^{M} \) is a one to one function
and \( |Dom(f_{i}^{M})|\leq \lambda _{1}[U] \).
\item For all \( i\leq m_{2} \), \( |s^{M}_{i}|\leq \lambda _{0}[U] \).
\end{enumerate}
\item We call \( \frak {M} \) a \noun{simple} model for \( \tau  \)
on \( \frak {U} \) if for all \( U\in \frak {U} \), \( \frak {M}[U] \)
is a simple model for \( \tau  \) on \( U \).
\item Let \( U\in \frak {U} \) and \( R \) a \( n(\frak {R}) \)-place
relation on \( U \). We call \( M \) a \noun{simple} expansion
of \( R \) on \( U \) for vocabulary \( \tau  \) if: \noun{}

\begin{enumerate}
\item \( r\in \tau  \).
\item \( M \) is a model for \( \tau  \) on \( U \).
\item \( r^{M}=R \).
\item The restriction of \( M \) to \( \tau \setminus \{r\} \) is a simple
model for \( \tau \setminus \{r\} \) on \( U \). In particular \( \tau \setminus \{r\} \)
is a simple vocabulary.
\end{enumerate}
\item We call \( \frak {M} \) a \noun{simple} expansion of \( \frak {R} \)
(on \( \frak {U} \)) for vocabulary \( \tau  \), if for all \( U\in \frak {U} \),
\( \frak {M}[U] \) is a simple expansion of \( \frak {R}[U] \) for
\( \tau  \) on \( U \). 
\end{enumerate}
\begin{defn}
Let \( \tau =\{f_{0},...,f_{m_{1}},s_{0},...,s_{m_{2}},c_{0},...,c_{m_{3}}\} \)
be a vocabulary, and \( \Delta  \) a set of formulas in \( \tau  \).
Let \( M \) be a model for \( \tau  \) on \( U \), \( m\in \omega  \),
\( A\subseteq U \), and \( \bar{a}\in {}^{m}U \). Define:
\begin{enumerate}
\item The \( \Delta  \)-type of \( \bar{a} \) over \( A \) in \( M \)
is: \[
tp_{\Delta }(\bar{a},A,M):=\{\varphi (\bar{x},\bar{b}):\varphi (\bar{x},\bar{y})\in \Delta ,lg(\bar{x})=m,\bar{b}\in {}^{<\omega }A,M\models \varphi (\bar{a},\bar{b})\}\]
 
\item \( S_{\Delta }^{m}(A,M):=\{tp_{\Delta }(\bar{a},A,M):\bar{a}\in {}^{m}U\} \).
if \( M=(U,R) \) we write \( S_{\Delta }^{m}(A,R) \) instead and
similarly in 1.
\item If \( p\in S_{\Delta }^{m}(A,M) \), \( \overline{a'}\in {}^{m}U \)
and \( \varphi (\bar{x},\bar{b})\in p\Rightarrow M\models \varphi (\overline{a'},\bar{b}) \),
then we say that \( \overline{a'} \) realizes \( p \). in particular
\( \bar{a} \) realizes \( tp_{\Delta }(\bar{a},A,M) \).
\end{enumerate}
\begin{defn}
\label{E_A} Let \( \tau =\{f_{0},...,f_{m_{1}},s_{0},...,s_{m_{2}},c_{0},...,c_{m_{3}}\} \)
be a vocabulary and \( \Delta  \) a set of formulas in \( \tau  \).
\begin{enumerate}
\item For all \( U\in \frak {U} \), \( A\subseteq U \) and \( M \) a
model for \( \tau  \) on \( U \), define an equivalence relation
\( E=E^{\Delta ,M}_{A,U} \) (we usually write \( E^{\Delta ,M}_{A} \)
where \( U \) is understood) on \( U \) by:\[
E:=\{(x',x'')\in {}^{2}U:tp_{\Delta }(x',A,M)=tp_{\Delta }(x'',A,M)\}\]

\item Let \( U\in \frak {U} \), \( m\in \omega  \) and \( E \) an equivalence
relation on \( U \). We call \( E \) \( m \)-big, if \( E \) has
at least \( m \) equivalence classes of size at least \( m \). If
\( E \) is not \( m \)-big we say it is \( m \)-small.
\item Let \( \frak {M} \) be a model for \( \tau  \) over on \( \frak {U} \),
define a function from \( \frak {U} \) to \( \mathbb {N} \), \( k_{\Delta }=k_{\Delta ,\frak {M}} \)
as follows: \( k_{\Delta }[U] \) is the maximal number \( k \) such
that there exists \( A\subseteq U \), \( |A|\leq \lambda _{0}[U] \),
and \( E_{A}^{\Delta ,\frak {M}[U]} \) is \( k \)-big. 
\end{enumerate}
\end{defn}
\end{defn}
\end{defn}
\begin{lem}
\label{int E}Let \( \frak {M} \) be a simple expansion of \( \frak {R} \)
for a vocabulary \( \tau  \), and \( \Delta  \) a finite set of
formulas in \( \tau  \). then: \( \{Q_{^{_{\frak {R}},}}Q_{\lambda _{0}}^{mon},Q_{\lambda _{1}}^{1-1}\}\geq _{int}Q^{eq}_{k_{\Delta ,\frak {M}},k_{\Delta ,\frak {M}}} \).
\end{lem}
\begin{proof}
For all \( U\in \frak {U} \), let \( A_{U}\subseteq U \) be the
subset the existence of which is promised by \ref{E_A}.3. Let \( s' \)
be a \( 1 \)-place relation symbol. Define a simple vocabulary \( \tau ':=\tau \cup \{s'\} \),
and a formula in \( \tau ' \):\[
\psi (x',x''):=(\forall \overline{b})\bigwedge _{\varphi (x,\overline{y})\in \Delta }\{s'(\overline{b})\rightarrow [\varphi (x',\overline{b})\equiv \varphi (x'',\overline{b})]\}\]
 (where \( (\forall \bar{b}) \) stands for \( \forall b_{0}...\forall b_{lg(\bar{b})-1} \),
and \( s'(\bar{b}) \) stands for \( \bigwedge _{i<lg(\bar{b})-1}s'(b_{i}) \)).
Let \( \frak {M}' \) be the simple expansion of \( \frak {M} \)
for \( \tau ' \) defined by \( s'^{\frak {M}'[U]}:=A_{U} \), for
all \( U\in \frak {U} \). Then for all \( U\in \frak {U} \) and
\( a,b\in U \): \[
aE^{\Delta ,\frak {M}[U]}_{A_{U}}b\Longleftrightarrow \frak {M}'[U]\models \psi (a,b)\]
 Define \( \frak {E} \) by \( \frak {E}[U]=E_{A_{U},U}^{\Delta ,\frak {M}[U]} \).
Since \( \frak {M}' \) is a simple expansion of \( \frak {R} \)
we have \( Q_{\frak {E}}\leq _{int}\{Q_{^{_{\frak {R}},}}Q_{\lambda _{0}}^{mon},Q_{\lambda _{1}}^{1-1}\} \)
(see \ref{int exp}) when the interpreting formula is \( \psi  \).
Now by \ref{Int k,k} we have \( Q_{\frak {E}}\geq _{int}Q^{eq}_{k_{\Delta ,\frak {M}},k_{\Delta ,\frak {M}}} \),
so by transitivity of \( \leq _{int} \) we are done.
\end{proof}
\begin{lem}
Let \( n \) be a natural number no larger than \( n(\frak {R}) \).
Let \( \tau  \) be a simple vocabulary, and \( \frak {M} \) a simple
expansion of \( \frak {R} \) for \( \tau \cup \{r\} \). Let \( \Delta  \)
be a finite set of formulas in \( \tau \cup \{r\} \), of the form
\( \varphi (x,\bar{y}) \) such that \( lg(\overline{y})\leq n \).
Let \( U\in \frak {U} \) and \( k=k_{\Delta ,\frak {M}}[U] \). Then
there exists \( A\subseteq U \) such that: 
\begin{enumerate}
\item \( |A|\leq n(k+1) \).
\item If \( \varphi (x,\overline{y})\in \Delta  \) and \( \overline{a}\in {}^{lg(\overline{y})}U \)
are a formula and parameters, then the formula \( \varphi (-,\overline{a}) \)
divides every equivalence class of \( E^{\Delta ,\frak {M}[U]}_{A} \)
into two parts one of which has no more than \( (k+1)\cdot 2^{m^{*}} \)
elements, where \( m^{*}=|\Delta |(k+1)^{n+1}\cdot n^{n} \).
\item There exists at least \( k+1 \) types, \( p\in S^{1}_{\Delta }(A,\frak {M}[U]) \),
realized by at least \( k\cdot 2^{m^{*}} \) elements of \( U \)
each.
\end{enumerate}
\end{lem}
\begin{proof}
Define a natural number \( m_{l} \) by dawnword induction on \( l\leq k+1 \):
\( m_{k+1}=0 \), \( m_{l}=|\Delta |(n(l+1))^{n}+m_{l+1} \). By induction
on \( l\leq k+1 \) we try to build a set \( A_{l}\subseteq U \)
such that \( |A_{l}|\leq n*l \), and there exists at least \( l \)
types \( p\in S^{1}_{\Delta }(A_{l},\frak {M}[U]) \) realized by
at least \( (k+1)*2^{m_{l}} \) elements each. If we succeed then
the existence of \( A_{k+1} \) is a contradiction to the definition
of \( k \). (We assume here that \( |A_{k+1}|\leq \lambda _{0}[U] \),
but without loss of generality we can assume that as \( |A_{k+1}| \)
is bounded. see also remark \ref{chenge lambda}). Let \( l_{0}<k+1 \)
be such that we have built \( A_{0},...,A_{l_{0}} \) but we can not
build \( A_{l_{0}+1} \). Put \( A=A_{l_{0}} \). Clearly \( A \)
satisfies (1). We prove (2).

Put \( M:=\frak {M}[U] \). Let \( \left\langle B_{i}:i\leq l_{0}\right\rangle  \)
be am enumeration of equivalence classes of \( E^{\Delta ,M}_{A} \)
with at least \( (k+1)*2^{m_{l_{0}}} \) elements. (note that there
are exactly \( l_{0} \) such classes since \( l_{0} \) is maximal).
Let \( \varphi (x,\overline{y})\in \Delta  \) and \( \textrm{ }\overline{a}\in {}^{lg(\overline{y})}U\textrm{ } \)
be some formula and parameters. The relation \( E^{\Delta }_{A\cup \overline{a}} \)
divides every class \( B_{i} \) to at most \( 2^{|\Delta |*(|A|+n)^{n}} \)
parts. Hence by the pigeon hole principle at least one of those parts
has at least \( \frac{|B_{i}|}{2^{|\Delta |*(|A|+n)^{n}}}\geq  \)\( \frac{|B_{i}|}{2^{|\Delta |*(n(l_{0}+1))^{n}}} \)\( \geq (k+1)*2^{m_{l_{0}+1}} \)
elements. If for some \( i \) there are more than one part with more
than \( (k+1)*2^{m_{l_{0}+1}} \) elements then define \( A_{l_{0}+1}=A\cup \overline{a} \)
and we get: 
\begin{enumerate}
\item \( |A_{l_{0}+1}|\leq |A_{l_{0}}|+|\overline{a}|\leq n*l_{0}+n\leq n(l_{0}+1) \).
\item There exists at least \( l_{0}+1 \) types \( p\in S^{1}_{\Delta }(A_{l_{0}+1},\frak {M}[U]) \)
realized by at least \( (k+1)\cdot 2^{m_{l_{0}+1}} \) elements each.
\end{enumerate}
This is a contradiction to the maximallity of \( l_{0} \). Now assume
towards contradiction that \( \varphi (-,\overline{a}) \) divides
some \( B_{i} \) into two parts, both larger than \( (k+1)*2^{m^{*}} \)
(note that \( m^{*}\geq m_{l_{0}} \) so there is no need to check
classes smaller than \( (k+1)*2^{m^{*}} \)). Then \( E_{A\cup \overline{a}} \)
divides each part into at most \( 2^{|\Delta |*(n(l_{0}+1))^{n}} \)
classes and hence each part contains an equivalence class of \( E_{A\cup \overline{a}} \)
with at least \( \frac{(k+1)*2^{m^{*}}}{2^{|\Delta |*(n(l_{0}+1))^{n}}}\geq  \)
\( \frac{(k+1)*2^{m_{l_{0}}}}{2^{|\Delta |*(n(l_{0}+1))^{n}}}= \)
\( (k+1)*2^{m_{l_{0}+1}} \) elements, so \( B_{i} \) contains two
such classes and this, as we saw, is a contradiction. To prove (3)
we note that \( l_{0}+1\leq k+1 \), and \( m^{*}\geq m_{l_{0}}\geq m_{l_{0}+1} \)
hence the existence of \( k+1 \) classes with \( k*2^{m^{*}} \)
elements contradicts the maximality of \( l_{0} \).
\end{proof}
\begin{thm}
\label{main}Let \( \tau  \) be a simple vocabulary, and \( \Delta  \)
a finite set of formulas in \( \tau \cup \{r\} \). \underbar{Then}
one of the following conditions hold:
\begin{enumerate}
\item There exists a sequence of worlds: \( \left\langle U_{i}\in \frak {U}:i\in \omega \right\rangle  \),
and a sequence of natural numbers: \( \left\langle n_{i}:i\in \omega \right\rangle  \)
such that \( n_{i}\longrightarrow \infty  \) and there exists a simple
vocabulary \( \tau ' \), a formula \( \varphi (x,y) \) in \( \tau '\cup \{r\} \)
and a simple expansion \( \frak {M}' \) of \( \frak {R} \) for \( \tau '\cup \{r\} \),
such that for all \( i\in \omega  \): \( \{(x,y)\in {}^{2}U_{i}:\frak {M}'[U_{i}]\models \varphi (x,y)\} \)
is an \( n_{i} \)-big equivalence relation on \( U_{i} \). 
\item There exists a natural number \( k^{*} \) such that for all \( U\in \frak {U} \)
and \( M \) - a simple expansion of \( \frak {R}[U] \) for \( \tau \cup \{r\} \)
on \( U \), there exist \( A=A^{\Delta ,M}_{U}\subseteq U \) such
that \( |A|\leq k^{*} \), \( E^{\Delta ,M}_{A,U} \) is \( k^{*} \)-small,
and for every formula \( \varphi (x,\overline{y})\in \Delta  \) and
parameters \( \overline{a}\in {}^{lg(\bar{y})}U \), \( \varphi (-,\overline{a}) \)
divides each equivalence class of \( E^{\Delta ,M}_{A} \) into two
parts one of which has less than \( k^{*} \) elements.
\end{enumerate}
\end{thm}
\begin{proof}
Define \( \mathbb {M} \) to be the class of all simple expansions
of \( \frak {R} \) for \( \tau \cup \{r\} \) on \( \frak {U} \).
For all \( U\in \frak {U} \) define:\[
k^{max}_{\Delta }[U]=max\{k_{\Delta ,\frak {M}}[U]:\frak {M}\in \mathbb {M}\}\]
 Note that the definition of \( k_{\Delta }^{max}[U] \) depends only
on \( \frak {M}[U] \) and not the over values of \( \frak {M} \)
(that is for each \( \frak {M}\in \mathbb {M} \)) and since \( |\{\frak {M}[U]:\frak {M}\in \mathbb {M}\}|<\aleph _{0} \)
the maximum is obtained. Next we assume that \( sup\{k^{max}_{\Delta }[U]:U\in \frak {U}\}=\aleph _{0} \)
and show that condition (1) is satisfied. Let \( \left\langle U_{i}\in \frak {U}:i\in \omega \right\rangle  \)
be a sequence of universes such that \( k^{max}_{\Delta }[U_{i}]\longrightarrow \infty  \),
and for all \( i\in \omega  \) define \( n_{i}=k^{max}_{\Delta }[U_{i}] \).
Define a simple vocabulary \( \tau '=\tau \cup \{s'\} \) (\( s' \)
a \( 1 \)-place relation symbol). We now define \( \frak {M} \)\( . \)
For all \( i\in \omega  \) note by \( M_{i} \) the model for \( \tau \cup \{r\} \)
on \( U \), for which the maximum in the definition of \( k^{max}_{\Delta }[U_{i}] \)
is obtained. Define \( \frak {M}[U_{i}]|\tau \cup \{r\}:=M_{i} \).
By \ref{E_A}.3 there exists a subset \( A_{i}\subseteq U_{i} \)
such that \( E^{\Delta ,M_{i}}_{A_{i}} \) is a \( n_{i} \)-big equivalence
relation on \( U \). Let \( s'^{\frak {M}[U_{i}]}=A_{i} \). That
defines \( \frak {M} \) (obviously the definition on universes not
among the \( U_{i} \) is irrelevant). We define \( \varphi (x,y) \)
to be the formula interpreting \( E^{\Delta ,M_{i}}_{A_{i}} \) (see
\ref{int E}) namely:\[
\varphi (x,y):=(\forall \overline{b})\bigwedge _{\psi (x,\overline{z})\in \Delta }\{s'(\overline{b})\rightarrow [\psi (x,\overline{b})\equiv \psi (y,\overline{b})]\}\]
 it is clear that condition (1) is satisfied.

We now assume that \( \{k^{max}_{\Delta }[U]:U\in \frak {U}\} \)
is bounded and let \( k \) by its bound. we show that condition (2)
is satisfied. Let \( n:=max\{lg(\bar{y}):\varphi (x,\bar{y})\in \Delta \} \).
We define \( k^{*}=Max\{(k+1)*2^{|\Delta |(k+1)^{n+1}*n^{n}},n(k+1)\} \).
Now let \( U\in \frak {U} \), and \( M \) a simple expansion of
\( \frak {R}[U] \) on \( U \) for vocabulary \( \tau \cup \{r\} \).
Let \( A\subseteq U \) be the subset the existence of which is promised
by the previous lemma. Then all the demands of (2) are clear from
the previous claim and the fact \( k\geq k_{\Delta }^{max}[U]\geq k_{\Delta ,\frak {M}}[U] \).
\end{proof}

\section{The Complicated Case of The Dichotomy}

In this section we assume that \( \frak {U} \) and \( \frak {R} \)
satisfy condition (1) in \ref{main}, that is we can uniformly interpret
an arbitrarily large equivalence relation. We show that in this case
we can interpret bounded number theory in the logic \( L(Q_{\frak {R}}) \).
It follows that the set of logicly valid sentences in \( L(Q_{\frak {R}}) \)
is not recursive. 

The following result is well known:

\begin{lem}
\label{3.1}Let \( E \) be an \( n \)-big equivalence relation on
a universe \( U \).  Then we can uniformly (that is using formulas
independent of \( U \) and \( E \)) interpret the model \( \left( \{0,1,...,n-1\};0,S,+,*\right)  \)
using a finite number of isomorphic copies of \( E \).
\end{lem}
\begin{cor}
In theorem \ref{main} if condition (1) is satisfied then we can uniformly
interpret number theory bounded by \( n_{i} \) using a finite number
of isomorphic copies of \( R[U_{i}] \).
\end{cor}
\begin{proof}
Straight from \ref{3.1}.
\end{proof}

\section{The Simple Case of The Dichotomy}

In this section we will interpret \( Q_{\frak {R}} \) when \( \frak {R} \)
is {}``simple{}'' that is when condition (1) in theorem \ref{main}
is not satisfied. We will show that in this case there exists a simple
model on \( \frak {U} \) in which it is possible to interpret \( \frak {R} \)
by a first order formula. In fact we prove \( Q_{\frak {R}}\leq _{int}\{Q_{\lambda _{0}}^{mon},Q_{\lambda _{1}}^{1-1}\} \)
so we get a full understanding of \( Q_{\frak {R}} \).

\subsection{Formalizing The Assumptions And The Main Theorem}

\begin{asmp}
\label{assumption}In this section we assume that \( \frak {U} \)
and \( \frak {R} \) do not satisfy condition (1) in theorem \ref{main}.
(Note that this condition is independent of \( \Delta ). \) Hence
from that theorem we get the following:
\begin{enumerate}
\item For every simple vocabulary \( \tau  \), and \( \Delta  \) a finite
set of formulas in \( \tau \cup \{r\} \), there exists a number \( k_{1}^{*}=k_{1}^{*}(\Delta ) \)
and a function that assigns to every \( U\in \frak {U} \) and \( M \)
- a simple expansion of \( \frak {R}[U] \) for \( \tau \cup \{r\} \)
on \( U \), a set \( A=A^{\Delta ,M}_{U}\subseteq U \) such that
condition (2) in theorem \ref{main} is satisfied, that is: \begin{itemize}

\item[(*)] \( |A|\leq k_{1}^{*} \), \( E^{\Delta ,M}_{A} \) is \( k_{1}^{*} \)-small,
and for every formula \( \varphi (x,\overline{y})\in \Delta  \) and
parameters \( \overline{a}\in {}^{lg(\bar{y})}U \), \( \varphi (-,\overline{a}) \)
divides each equivalence class of \( E^{\Delta ,M}_{A} \) into two
parts one of which has at most \( k_{1}^{*} \) elements. \end{itemize} 
\item For every simple vocabulary \( \tau  \), and every formula \( \varphi (x,y) \)
in \( \tau \cup \{r\} \), there exists a natural number \( k_{2}^{*}=k_{2}^{*}(\varphi ) \)
such that: \begin{itemize}
\item[(**)] If \( \frak {M} \) is a simple expansion of \( \frak {R} \)
for \( \tau \cup \{r\} \) and \( U\in \frak {U} \), then the interpretation
of \( \varphi (x,y) \) in \( \frak {M}[U] \) (that is \( \{(x,y)\in {}^{2}U:\frak {M}[U]\models \varphi (x,y)\} \))
is \noun{not} a \( k^{*}_{2} \)-big equivalence relation. \end{itemize}
\end{enumerate}
\end{asmp}
\begin{rem}
\label{inc k_1}We can increase \( k^{*}_{1}(\Delta ) \), meaning
if \( m\geq k^{*}_{1}(\Delta ) \) than \( m \) satisfies \( (*) \)
(for the same function \( A_{U}^{\Delta ,M} \)). Hence: 
\begin{enumerate}
\item If we are given a function \( \Delta \mapsto m(\Delta ) \) then without
loss of generality (by changing the definition of \( k^{*}_{1} \))
we may assume: \( k^{*}_{1}(\Delta )\geq m(\Delta ) \) for all \( \Delta  \).
\item If \( \Delta \subseteq \Delta ' \) then without loss of generality
(by redefining \( k^{*}_{1} \) by induction on \( |\Delta | \))
we may assume \( k^{*}_{1}(\Delta ')\geq k^{*}_{1}(\Delta ) \).
\end{enumerate}
\end{rem}
We now formalize the main theorem of this section.

\begin{thm}
\label{simple} There exists a simple vocabulary \( \tau  \), and
a first order formula \( \varphi (x_{0},...,x_{n(\frak {R})-1}) \)
in \( \tau  \), and there exists \( \frak {M} \) a simple expansion
of \( \frak {R} \) for \( \tau \cup \{r\} \) on \( \frak {U} \)
such that for all \( U\in \frak {U} \):\[
\frak {M}[U]\models (\forall \bar{x})[r(\bar{x})\equiv \varphi (\bar{x})]\]
  
\end{thm}
\begin{cor}
\( Q_{\frak {R}}\leq _{int}\{Q^{1-1}_{\lambda _{1}},Q^{mon}_{\lambda _{0}}\} \).
\end{cor}
\begin{proof}
Straight from the theorem when the interpreting formula is \( \varphi  \).
\end{proof}
In the rest of the paper we will prove theorem \ref{simple}.

\subsection{Proof of The Main Theorem in The \protect\( 2\protect \)-place Case}

We prove theorem \ref{simple} under the assumption \( n(\frak {R})=2 \).
\( \Delta  \) will be a finite set of formulas with at most \( 2 \)
free variables in the vocabulary \( \{r\} \). In other words \( \tau =\emptyset  \).
Hence the set \( A_{U}^{\Delta ,M} \) (see \ref{assumption}) and
the relation \( E_{A}^{\Delta ,M} \) (see \ref{E_A}) are independent
of \( M \), and depend on \( \Delta  \) alone so they will by noted
by \( A^{\Delta }_{U} \) and \( E_{A}^{\Delta } \). 

\begin{defn}
\label{S}Let \( \Delta  \) be as above and \( U\in \frak {U} \).
Let \( k^{*}=k^{*}_{1}(\Delta ) \) and \( A=A^{\Delta }_{U} \) be
the ones we get from \ref{assumption}.1. We define:
\begin{enumerate}
\item For all \( \varphi (x,y)\in \Delta  \) and \( y_{0}\in U \): \[
Minority_{\Delta }(y_{0},\varphi ):=\{x_{0}\in U:|\{x\in U:xE^{\Delta }_{A}x_{0}\wedge \varphi (x,y_{0})\equiv \varphi (x_{0},y_{0})\}|\leq k^{*}\}\]

\item \( S=S^{\Delta } \) is the \( 2 \)-place relation on \( U \) given
by:\[
x_{0}Sy_{0}\Leftrightarrow x_{0}\in \bigcup _{\varphi (x,y)\in \Delta }Minority_{\Delta }(y_{0},\varphi )\]

\end{enumerate}
\end{defn}
\begin{lem}
\label{bound}Let \( \Delta  \) be as above. We use the notations
of the previous definition and also note \( k_{2}^{*}=k_{2}^{*}(\psi ) \)
(see \ref{assumption}.2) where \( \psi (x',x''):=(\forall \overline{b})\bigwedge _{\varphi (x,\overline{y})\in \Delta }\{s(\overline{b})\rightarrow [\varphi (x',\overline{b})\equiv \varphi (x'',\overline{b})]\} \).
Then:
\begin{enumerate}
\item \( \left| \left\{ x:\left| x/E^{\Delta }_{A}\right| \leq 2\cdot k^{*}\right\} \right| \leq k^{*}\cdot 2^{|\Delta |k^{*}+1} \),
and we write \( l^{*}=k^{*}\cdot 2^{|\Delta |k^{*}+1} \).
\item For all \( y\in U \): \( \left| \left\{ x:xSy\right\} \right| \leq |\Delta |\cdot \left( k^{*}\right) ^{2}+l^{*} \).
\item \( \left| \left\{ x:\left| \left\{ y:xSy\right\} \right| >2^{|\Delta |(k^{*}+k_{2}^{*})}\cdot k_{2}^{*}+l^{*}\right\} \right| \leq |\Delta |\cdot (k_{2}^{*}\cdot k^{*})^{2}\cdot 2^{|\Delta |(k^{*}+k^{*}_{2})}+l^{*} \).
\end{enumerate}
\end{lem}
\begin{proof}
(1): The number of types \( p\in S^{1}_{\Delta }(A,\frak {R}[U]) \)
is no larger than \( 2^{|\Delta ||A|} \) since for every formula
in \( \Delta  \) there are at most two free variables. We also have
\( |A|\leq k^{*} \). So the number of equivalence classes of \( E_{A}^{\Delta } \)
is no larger than \( 2^{|\Delta |k^{*}} \) and (1) follows directly.

(2): Let \( x,y\in U \). Assume \( |x/E^{\Delta }_{A}|>2\cdot k^{*} \).
For all \( \varphi \in \Delta  \) we have \( x/E^{\Delta }_{A}\cap Minority_{\Delta }(y,\varphi )\leq k^{*} \).
Hence \( |\{x':xE^{\Delta }_{A}x'\wedge x'S^{\Delta }y\}|\leq |\Delta |\cdot k^{*} \).
The number of equivalence classes of \( E^{\Delta }_{A} \) which
are larger than \( 2\cdot k^{*} \) is also no larger than \( k^{*} \).
Hence we get: \( |\{x:|x/E^{\Delta }_{A}|>2\cdot k^{*}\wedge xS^{\Delta }y\}|\leq |\Delta |\cdot (k^{*})^{2} \).
To this we add at most \( l^{*} \) elements from {}``small classes{}''
(see (1)) and (2) follows.

(3): We write \( m=|\Delta |\cdot (k_{2}^{*}\cdot k^{*})^{2}\cdot 2^{|\Delta |(k^{*}+k^{*}_{2})} \).
First we disregard all the elements of \( \left\{ x:\left| x/E_{A}\right| \leq 2\cdot k^{*}\right\}  \)
and using (1) we decrease the bounds by \( l^{*} \). So seeking a
contradiction we assume that there are different \( \{x_{0},...,x_{m}\} \)
so that for each \( i\leq m \) there exists different \( \{y^{i}_{0},...,y_{2^{|\Delta |(k^{*}+k_{2}^{*})}\cdot k_{2}^{*}}^{i}\} \)
with \( x_{i}Sy^{i}_{j} \). Using (2) (with the bounds in (2) also
decreased by \( l^{*} \)) there exists a subset of \( \{x_{0},...,x_{m}\} \)
with at least \( k^{*}_{2} \) elements such that the elements of
\( Y_{x_{i}}:=\{y^{i}_{j}:x_{i}Sy^{i}_{j}\} \) are pairwise disjoint
(see figure). Without loss of generality we assume that this set is
\( \{x_{0},...,x_{k_{2}^{*}}\} \). For every \( x_{i} \) the sets
of \( Y_{x_{i}} \) satisfy at most \( 2^{|\Delta |(k^{*}+k^{*}_{2})} \)
different types \( p\in S^{1}_{\Delta }(A\cup \{x_{i}:i\leq k^{*}_{2}\},\frak {R}[U]) \).
Hence there are more than \( k^{*}_{2} \) of them that satisfy the
same type (again we assume those are the first elements). In conclusion
we get \( \{x_{i}\}_{i=0}^{k_{2}^{*}} \) and \( \{y^{i}_{j}\}_{i,j=0,...,k^{*}_{2}} \)
without repetitions such that the type \( tp_{\Delta }(y^{i}_{j},A\cup \{x_{i}:i\leq k^{*}_{2}\},\frak {R}[U]) \)
is independent of \( j \), and \( x_{i_{1}}Sy_{j}^{i_{2}}\Leftrightarrow i_{1}=i_{2} \)
holds for all \( i_{1},i_{2},j\leq k^{*}_{2} \). So \( \psi  \)
with \( s \) taken to represent \( A\cup \{x_{i}:i\leq k^{*}_{2}\} \)
interprets a \( k^{*}_{2} \)-big equivalence relation on \( U \).
This is a contradiction to the definition of \( k^{*}_{2} \).
\end{proof}
\vspace{0.3cm}
{\centering \resizebox*{0.9\textwidth}{0.2\textheight}{\includegraphics{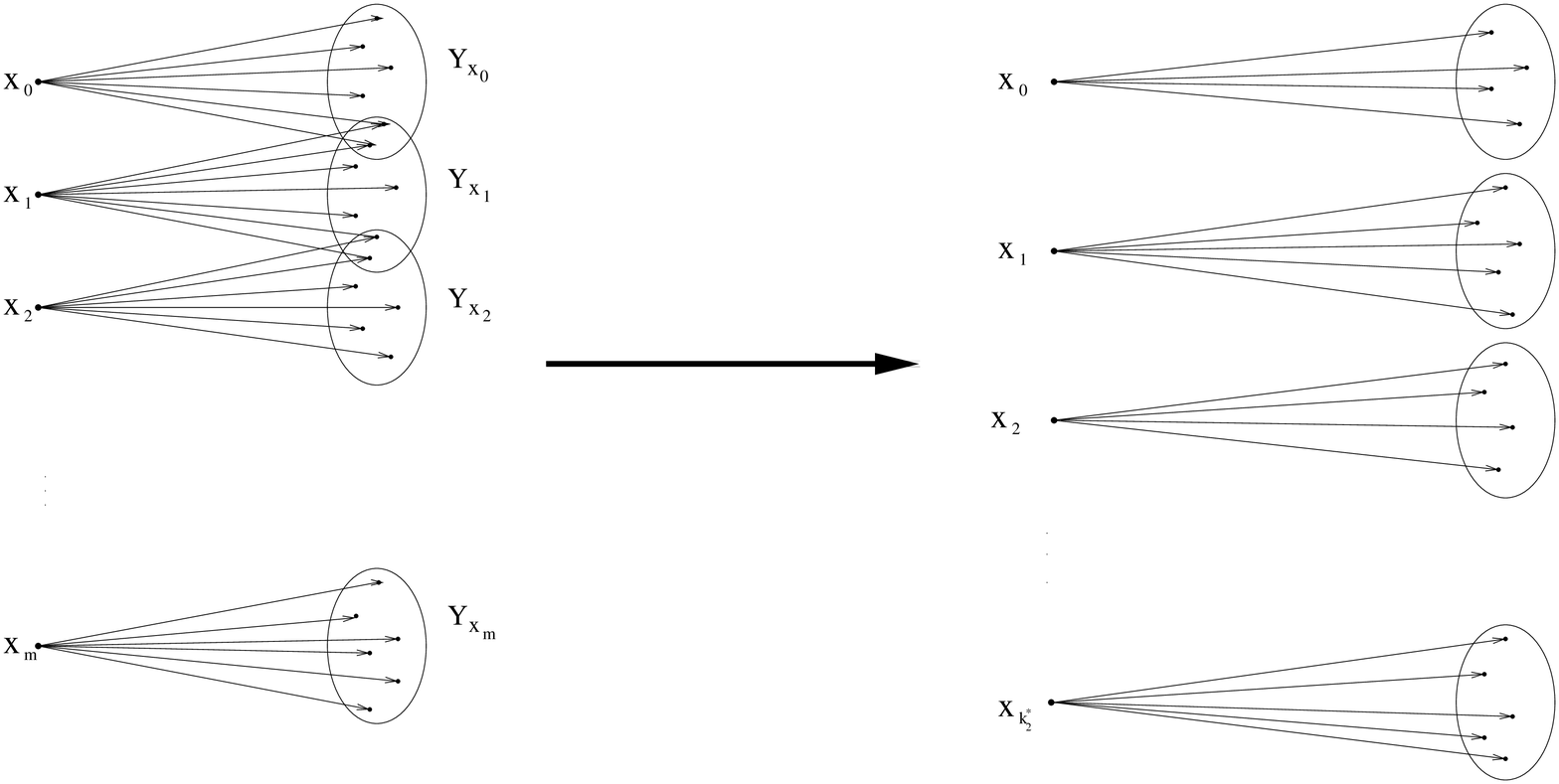}} \par}
\vspace{0.3cm}

\begin{lem}
\label{int R n=3D2}There exist a simple vocabulary \( \tau  \),
and a finite set of formulas \( \Phi  \) in \( \tau  \), and a simple
model \( \frak {M} \) for \( \tau  \) on \( \frak {U} \), such
that for all \( U\in \frak {U} \) and \( x,x',y,y'\in U \) if \( tp_{\Phi }((x,y),\emptyset ,\frak {M}[U])=tp_{\Phi }((x',y'),\emptyset ,\frak {M}[U]) \)
then \( (U,\frak {R}[U])\models r(x,y)\equiv r(x',y') \). 
\end{lem}
\begin{proof}
We simultanios define \( \tau  \) and its interpretation \( \frak {M}[U] \)
for some \( U\in \frak {U} \). \( \Phi  \) will be the set of atomic
formulas in \( \tau  \) with terms of the form \( x,f(x),c,f(c) \)
(function composition is not allowed). for gravity we write: \( M:=\frak {M}[U] \)
and \( R:=\frak {R}[U] \). Let \( \Delta :=\{r(x,y)\} \). Using
the notations of \ref{bound} we define:\[
A^{*}=A\cup \{x:|x/E^{\Delta }_{A}|\leq 2\cdot k^{*}\}\cup \left\{ x:\left| \left\{ y:xS^{\Delta }y\right\} \right| >2^{|\Delta |(k^{*}+k_{2}^{*})}\cdot k_{2}^{*}+l^{*}\right\} \]
 by \ref{bound} \( |A^{*}| \) is uniformly bounded (that is the
bound is independent of \( U \)). \( \tau  \) will contain: private
constants for all the elements of \( A^{*} \)(\( \{c_{x}:x\in A^{*}\} \),
\( c_{x}^{M}:=x \)), and \( 1 \)-place relation symbols for the
equivalence classes of \( E^{\Delta }_{A^{*}} \) (\( \{s_{x/E^{\Delta }_{A^{*}}}:x\in U\} \),
\( s^{M}_{x/E^{\Delta }_{A^{*}}}:=x/E^{\Delta }_{A^{*}} \)). Note
that the number of such classes is also uniformly bounded. Now We
look at \( S^{\Delta }|U\setminus A^{*} \) this is a digraph with
(uniformly) bounded degree, that is for all \( x\in U\setminus A^{*} \),
\( |\{y\notin A^{*}:xS^{\Delta }y\}| \) is bounded by \( 2^{|\Delta |(k^{*}+k_{2}^{*})}\cdot k_{2}^{*}+l^{*} \)
and for all \( y\in U\setminus A^{*} \), \( |\{x\notin A^{*}:xS^{\Delta }y\}| \)
is bounded by \( |\Delta |*\left( k^{*}\right) ^{2}+l^{*} \) (see
\ref{bound}). Hence we can divide \( S^{\Delta }|U\setminus A^{*} \)
into \( \left\langle S_{m}:m<m^{*}\right\rangle  \) with: \( \bigcup _{m<m^{*}}S_{m}=S^{\Delta }|U\setminus A^{*} \)
and for all \( m<m^{*} \), \( S_{m} \) is a digraph with degree
1, that is a one to one partial function on \( U\setminus A^{*} \).
Note that \( m^{*} \) is uniformly bounded, in fact it is bounded
by the sum of the two bounds mentioned above. We add to \( \tau  \),
\( 1 \)-place function symbols \( \{f_{m}:m<m^{*}\} \) and define
\( f_{m}^{M}:=S_{m} \). 

Let \( \left\langle B_{i}:i<i^{*}\right\rangle  \) be an enumeration
of \( \{x/E_{A}^{\Delta }\setminus A^{*}:|x/E^{\Delta }_{A}|>2\cdot k^{*}\} \).
Note that \( i^{*}\leq k^{*} \). For all \( y\in U \) and \( i<i^{*} \)
there is a truth value \( t^{y}_{i} \) that is the value the formula
\( r(-,y) \) gets for the \noun{majority} of the elements of \( B_{i} \).
Since we deal with {}``big{}'' classes (that is with more than \( 2\cdot k^{*} \)
elements) we get: for all \( y\in U \), \( i<i^{*} \) and \( x\in B_{i} \),
\( R(x,y)=t^{y}_{i}\Leftrightarrow \lnot xS^{\Delta }y \). We divide
each \( B_{i} \) into \( 2^{i^{*}} \) parts according to the truth
values, \( t^{y}_{i}:i<i^{*} \). This means that for each part, the
value of the vector \( \left\langle t^{y}_{i}:i<i^{*}\right\rangle  \)
is independent of \( y \). For all \( i<i^{*} \), we note these
parts by \( \left\langle B^{i}_{j}:j<2^{i^{*}}\right\rangle  \).
We add to \( \tau  \), \( 1 \)-place relations \( \{s_{i,j}:i<i^{*},j<2^{i^{*}}\} \)
and define \( s_{i,j}^{\frak {M}[U]}:=B^{i}_{j} \). This completes
the definition of \( \tau  \) and \( \frak {M} \). 

We now prove that \( \frak {M} \) is as desired. Let \( a,a',b,b'\in U \)
and assume \( tp_{\Phi }((a,b),\emptyset ,\frak {M}[U])=tp_{\Phi }((a',b'),\emptyset ,\frak {M}[U]) \).
If \( a\in A^{*} \) then \( a=a' \) (due to the formula \( x=c_{x} \)),
and the truth value of \( R(a,b) \) is determined by \( b/E_{A^{*}}^{\Delta } \).
Moreover \( b/E^{\Delta }_{A^{*}}=b'/E^{\Delta }_{A^{*}} \) (due
to the formula \( s_{b/E_{A^{*}}^{\Delta }}(y) \)), so we get \( R(a,b)=R(a',b') \)
as desired. Symmetricly we deal with the cases \( b,b'a'\in A^{*} \).
So we can assume \( a,a',b,b'\notin A^{*} \). By the definition of
the functions \( S_{m} \) we have:

\( aS^{\Delta }b\Longleftrightarrow (\exists m<m^{*})aS_{m}b \)

\( a'S^{\Delta }b'\Longleftrightarrow (\exists m<m^{*})a'S_{m}b' \)

But due to the formulas of the form \( f_{m}(x)=y \), the right hand
side of both equations is equivalent, so we have \( aS^{\Delta }b\Leftrightarrow a'S^{\Delta }b' \).
Assume \( a\in B^{i_{1}}_{j_{1}} \), \( b\in B^{i_{2}}_{j_{2}} \).
Due to the formula \( s_{i,j}(x) \) we get \( a'\in B^{i_{1}}_{j_{1}} \),
\( b'\in B^{i_{2}}_{j_{2}} \). By the construction of the \( B_{i}^{j} \)
we get: 

\( R(a,b)=t^{b}_{i_{1}}\Longleftrightarrow \lnot aS^{\Delta }b \)

\( R(a',b')=t^{b'}_{i_{1}}\Longleftrightarrow \lnot a'S^{\Delta }b' \)

But \( b,b'\in B^{i_{2}}_{j_{2}} \) so \( t^{b}_{i_{1}}=t^{b'}_{i_{1}} \),
and as we have seen \( aS\Delta b\Leftrightarrow a'S^{\Delta }b' \).
Hence \( R(a,b)=R(a',b') \) as desired.
\end{proof}
\begin{cor}
\label{n=3D2} Theorem \ref{simple} is true for the case \( n(\frak {R})=2 \).
\end{cor}
\begin{proof}
Let \( \tau ,\Phi  \) and \( \frak {M} \) be as in the previous
claim. Note that by the definition of \( \Phi  \) all the formulas
in \( \Phi  \) have at most \( 2 \) free variables. Define \( C=\{\underline{t}|\underline{t}:\Phi \rightarrow \{\mathbb {T},\mathbb {F}\}\} \)
(each member of \( C \) represents a type in \( S_{\Phi }^{2}(\emptyset ,\frak {M}[U]) \)).
For all \( D\subseteq C \) define:\[
\chi _{D}(x,y):=\bigvee _{\underline{t}\in D}[\bigwedge _{\psi (x,y)\in \Phi }\psi (x,y)^{\underline{t}(\psi )}]\]
 where \( \psi ^{\mathbb {T}}=\psi  \), \( \psi ^{\mathbb {F}}=\neg \psi  \)
and \( \chi _{\emptyset }(x,y):=x\neq x \). (The formula \( \chi _{D}(x,y) \)
means \( (x,y) \) satisfies one of the types in \( D \)). For all
\( U\in \frak {U} \) define \( D_{U}\subseteq C \) by:\[
\{\underline{t}\in C|(\exists x,y\in U)[(U,\frak {R}[U])\models r(x,y)\wedge ((U,\frak {M}[U])\models \bigwedge _{\psi (x,y)\in \Phi }\psi (x,y)^{\underline{t}(\psi )})]\}\]
 This means \( D_{U} \) is the collection of types \( tp_{\Phi }((x,y),\emptyset ,\frak {M}[U]) \)
such that \( (U,\frak {R}[U])\models r(x,y) \). Using the previous
claim it is easy to verify that for all \( U\in \frak {U} \) and
\( x,y\in U \) we have:\[
(U,\frak {R}[U])\models r(x,y)\Longleftrightarrow (U,\frak {M}[U])\models \chi _{D_{U}}(x,y)\]
 We now add to \( \tau  \) constants: \( \{c_{true}\}\cup \{c_{D}:D\subseteq C) \).
For each \( U\in \frak {U} \), \( c_{true} \) is interpreted in
\( \frak {M}[U] \) by some element of \( U \). The rest of the constants
are interpreted so that for all \( D\subseteq C \): \( (c^{\frak {M}[U]}_{D}=c^{\frak {M}[U]}_{true})\Leftrightarrow (D=D_{U}) \)
holds. (assuming \( U \) has more than one element there is no problem
to do that). Now the desired formula in theorem \ref{simple} is:\[
\varphi (x,y):=\bigwedge _{D\subseteq C}[(c_{D}=c_{true})\rightarrow \chi _{D}(x,y)]\]

\end{proof}

\section{Proof of The Main Theorem in The General Case}

We prove theorem \ref{simple} when \( n(\frak {R})>2 \). from here
on we assume:

\begin{asmp}
\label{assumption2}\( \tau  \) is a simple vocabulary. \( \Delta  \)
is a finite set of formulas in \( \tau \cup \{r\} \), such that \( \varphi (\bar{x})\in \Delta \longrightarrow lg(\bar{x})\leq n(\frak {R}) \).
\end{asmp}
First we generalize definition \ref{S}.

\begin{defn}
Let \( \tau ,\Delta  \) be as above. Let \( U\in \frak {U} \) and
\( M \) be a simple expansion of \( \frak {R}[U] \) on \( U \)
for \( \tau \cup \{r\} \). Let \( n<n(\frak {R}) \). We note \( k^{*}=k_{1}^{*}(\Delta ) \)
and \( A=A^{\Delta ,M}_{U} \) the existence of which follows from
\ref{assumption}.1 and define:
\begin{enumerate}
\item For all \( \varphi (x,\overline{y})\in \Delta  \) with \( lg(\overline{y})=n \)
and \( \overline{b}\in {}^{n}U \):\[
Minority_{\Delta ,M}(\overline{b},\varphi ):=\{x\in U:|\{x'\in U:xE^{\Delta ,M}_{A}x'\wedge \varphi (x,\bar{b})\equiv \varphi (x',\bar{b})\}|\leq k^{*}\}\]

\item Define a relation \( S_{\Delta ,M}^{n}\subseteq U\times {}^{n}U \):\[
aS_{\Delta ,M}^{n}\overline{b}\Leftrightarrow a\in \bigcup \{Minority_{\Delta ,M}(\overline{b},\varphi ):\varphi (x,\overline{y})\in \Delta ,lg(\overline{y})=n\}\]

\end{enumerate}
\end{defn}
\begin{rem}
\label{s' imply s}For \( i\in \{1,2\} \) assume \( \tau _{i},\Delta _{i} \)
satisfy \ref{assumption2}, and \( \frak {M}_{i} \) is a simple expansion
of \( \frak {R} \) on \( \frak {U} \) for \( \tau _{i}\cup \{r\} \).
Furthermore assume \( \tau _{1}\subseteq \tau _{2} \), \( \Delta _{1}\subseteq \Delta _{2} \)
and \( \frak {M}_{1}=\frak {M}_{2}|\tau _{1} \). By \ref{inc k_1}.2
we may assume \( k^{*}_{1}(\Delta _{2})\geq k^{*}_{1}(\Delta _{1}) \),
hence for all \( U\in \frak {U} \) we can assume without loss of
generality (we can add elements to \( A^{\Delta ',\frak {M}'[U]}_{U} \)
if needed) that \( aS^{n}_{\Delta _{2},\frak {M}_{2}[U]}\bar{b}\Longrightarrow aS^{n}_{\Delta _{1},\frak {M}_{1}[U]}\bar{b} \). 
\end{rem}
\begin{lem}
Using the notations of the previous definition:
\begin{enumerate}
\item \( \left| \left\{ x:\left| x/E^{\Delta ,M}_{A}\right| \leq 2\cdot k^{*}\right\} \right| \leq k^{*}\cdot 2^{|\Delta |{k^{*}\choose n(\frak {R})}+1} \).
\item For all \( \overline{b}\in {}^{n}U \): \( \left| \left\{ x\in U:xS^{n}_{\Delta ,M}\overline{b}\right\} \right| \leq |\Delta |\cdot \left( k^{*}\right) ^{2}+k^{*}\cdot 2^{|\Delta |{k^{*}\choose n(\frak {R})}+1} \).
\item We write: \( l^{*}=l^{*}(\Delta ):=|\Delta |\cdot \left( k^{*}\right) ^{2}+k^{*}\cdot 2^{|\Delta |{k^{*}\choose n(\frak {R})}+1} \).
\end{enumerate}
\end{lem}
\begin{proof}
Similar to the proof of \ref{bound}, only in (1) we have at most
\( k^{*}\choose n(\frak {R}) \) different choices of parameters for
each formula.
\end{proof}
\begin{lem}
\textbf{\label{symetry}Symmetry Lemma (with Parameters):}

Assume \( \tau ,\Delta  \) satisfy \ref{assumption2}, and Let \( \frak {M} \)
be a simple expansion of \( \frak {R} \) for \( \tau \cup \{r\} \).
Let \( n<n(\frak {R}) \). Then there exists a simple vocabulary \( \tau '\supseteq \tau  \),
and \( \frak {M}' \) a simple expansion of \( \frak {M} \) for \( \tau '\cup \{r\} \),
and for \( i\in \{1,2\} \) there exists \( \Delta _{i}=\Delta _{i}(\Delta ) \)
such that \( \tau ',\Delta _{i} \) also satisfy \ref{assumption2}
and for all \( U\in \frak {U} \), \( a,b\in U \) and \( \bar{c}\in {}^{n-1}U \):\[
aS^{n}_{\Delta ,\frak {M}[U]}b\bar{c}\Longrightarrow (aS^{n-1}_{\Delta _{1},\frak {M}'[U]}\bar{c})\lor (bS^{n}_{\Delta _{2},\frak {M}'[U]}a\bar{c})\]

\end{lem}
\begin{proof}
First we define a few constants we will use later:

\( m^{*}:=k^{*}_{2}(\phi ) \) (see assumption \ref{assumption}).
where \( \phi  \) is the following formula in \( \tau \cup \{s,c_{1},...,c_{n-1}\}\cup \{r\} \):\[
\phi =\phi (y,y'):=(\forall x)\bigwedge _{\psi (x,y,z_{1},...,z_{n-1})\in \Delta }\{s(x)\rightarrow [\psi (x,y,c_{1},...,c_{n-1})\equiv \psi (x,y',c_{1},...,c_{n-1})]\}\]
 (\( s \) is a \( 1 \)-place relation symbol and \( c_{1},...,c_{n-1} \)
are constants not in \( \tau  \)). We also define: \( m_{1}=m_{1}(\Delta ):=(m^{*})^{2}\cdot 2^{|\Delta |m^{*}}\cdot l^{*}(\Delta ) \)
and \( m_{2}=m_{2}(\Delta ):=m^{*}\cdot 2^{|\Delta |m^{*}} \), where
\( l^{*}(\Delta ) \) was defined in the previous lemma.

let \( \frak {M}',\tau ' \) and \( \psi (x,x') \) be the vocabulary,
model and formula which interpret \( E_{A}^{\Delta ,M} \) (see the
proof of \ref{int E}). We define in \( \tau ' \) a formula that
will interpret \( xS^{n}_{\Delta ,\frak {M}[U]}y\bar{z} \) in \( \frak {M}'[U] \)
(where \( lg(\bar{z})=n-1 \)):\[
\chi (x,y,\bar{z}):=\bigvee _{\varphi (u,v,\bar{w})\in \Delta ,lg(\bar{w})=n-1}(\exists ^{\leq k^{*}_{1}(\Delta )}x')[\psi (x,x')\wedge (\varphi (x,y,\bar{z})\equiv \varphi (x',y,\bar{z}))]\]
 We therefor get: \begin{itemize}
\item[(*)] for all \( U\in \frak {U} \), \( a,b\in U \) and \( \bar{c}\in {}^{n-1}U \):
\( \frak {M}'[U]\models \chi (a,b,\bar{c})\Longleftrightarrow aS^{n}_{\Delta ,\frak {M}[U]}b\bar{c} \)

\end{itemize}

Define: \[
\chi '(x,\bar{z}):=(\exists ^{>m_{2}}y)\chi (x,y,\bar{z})\]
\[
\Delta _{1}:=\Delta \cup \{\chi '(x,\bar{z})\}\]
\[
\Delta _{2}:=\Delta \cup \{\chi (x,y,\bar{z})\}\]

Note that by \ref{inc k_1}.1 we may assume that \( k^{*}_{1}(\Delta )\geq max\{m_{1}(\Delta ),m_{2}(\Delta )\} \),
and by \ref{inc k_1}.2 we may assume \( k^{*}_{1}(\Delta _{i})\geq k^{*}_{1}(\Delta )\geq m_{i}(\Delta ) \)
for \( i\in \{1,2\} \). We now assume towards contradiction that
there exists \( U\in \frak {U} \), \( a,b\in U \) and \( \bar{c}\in {}^{n-1}U \)
such that:
\begin{enumerate}
\item \( aS^{n}_{\Delta ,\frak {M}[U]}b\bar{c} \).
\item \( \lnot (aS^{n-1}_{\Delta _{1},\frak {M}'[U]}\bar{c}) \).
\item \( \lnot (bS^{n}_{\Delta _{2},\frak {M}'[U]}a\bar{c}) \).
\end{enumerate}
\relax From (3) and \( k^{*}_{1}(\Delta _{2})\geq m_{2} \) we can find \( \{b_{0},...,b_{m_{2}}\} \)
without repetitions such that for all \( i\leq m_{2} \): \( \frak {M}'[U]\models \chi (a,b,\bar{c})\equiv \chi (a,b_{i},\bar{c}) \).
from (1) and \( (*) \) we get that for all \( i\leq m_{2} \): \( \frak {M}'[U]\models \chi (a,b_{i},\bar{c}) \).
Hence \( \frak {M}'[U]\models \chi '(a,\bar{c}) \). from (2) and
\( k^{*}_{1}(\Delta _{1})\geq m_{1} \) we can find \( \{a_{0},...,a_{m_{1}}\} \)
without repetitions such that for all \( i\leq m_{1} \): \( \frak {M}[U]\models \chi '(a,\bar{c})\equiv \chi '(a_{i},\bar{c}) \).
We have seen that \( \frak {M}'[U]\models \chi '(a,\bar{c}) \) so
by the definition of \( \chi '(x,\bar{z}) \) we have for all \( i\leq m_{1} \),
there exists \( \{b^{i}_{0},...,b^{i}_{m_{2}}\} \) without repetitions
such that \( i\leq m_{1}\land j\leq m_{2}\Rightarrow a_{i}S^{n}_{\Delta ,\frak {M}[U]}b^{i}_{j}\bar{c} \).
By the definition of \( l^{*}(\Delta ) \) and a repeated use of the
pigeon hole principle we can find a subset of \( \{a_{0},...,a_{m_{1}}\} \),
\( \{a_{i_{0}},...,a_{i_{m^{*}}}\} \) such that the sets \( \{\{b^{i_{l}}_{0},...,b^{i_{l}}_{m_{2}}\}:l\leq m^{*}\} \)
are pairwise disjoint. without loss of generality we assume \( i_{l}=l \)
for all \( l\leq m^{*} \). Using the pigeon hole principle again
we can find for all \( i\leq m^{*} \) subset of \( \{b^{i}_{0},...,b^{i}_{m_{2}}\} \)
with \( m^{*}+1 \) elements (and again we assume this subset is \( \{b^{i}_{0},...,b^{i}_{m^{*}}\} \))
such that for all \( \varphi (x,y,\bar{z})\in \Delta  \) and \( j_{1},j_{2}\leq m^{*} \)
we have \( \varphi (a_{i},b^{i}_{j_{1}},\bar{c})\Leftrightarrow \varphi (a_{i},b^{i}_{j_{2}},\bar{c}) \).
In conclusion we got: \( \{a_{0},...,a_{m^{*}}\} \) without repetitions
and for each \( i\leq m^{*} \): \( \{b^{i}_{0},...,b^{i}_{m^{*}}\} \)
without repetitions such that \( a_{i_{1}}S^{n}_{\Delta }b^{i_{2}}_{j}\bar{c}\Leftrightarrow i_{1}=i_{2} \).
Moreover the elements of \( \{b^{i}_{0},...,b^{i}_{m^{*}}\} \) satisfy
the same formulas of the form \( \varphi (a_{i},y,\bar{c})\in \Delta  \)
(\( \bar{c} \) and \( a_{i} \) are parameters). Now the formula
\( \phi (y,y') \) (where \( s \) is taken to mean \( \{a_{0},...,a_{m^{*}}\} \)
and the constants \( c_{i} \) are taken to mean the elements \( c_{i} \))
interprets a \( m^{*}+1 \)-big equivalence relation on \( \{y^{i}_{j}:i,j\leq m^{*}\} \).
This is a contradiction to the definition of \( m^{*} \). 
\end{proof}
We now prove a number of lemmas we need for the proof of the main
theorem. First we show that we can code a delta system of \( n \)-tuples
by singletons:

\begin{lem}
\label{delta}Let n be a natural number. Then there exists a simple
vocabulary \( \tau  \), and a formula \( \theta (x,\bar{y}) \) in
\( \tau  \) with \( lg(\bar{y})=n \) such that: for all \( U\in \frak {U} \)
and delta system \( \left\langle \overline{a^{i}}\in {}^{n}U:i<i^{*}\right\rangle  \)
(\( i^{*} \) some natural number), we have a simple model \( M \)
for \( \tau  \) on \( U \) and a sequence \( \left\langle b_{i}\in U:i<i^{*}\right\rangle  \)
such that:\[
(\forall \bar{a}\in {}^{n}U)(\forall b\in U)[M\models \theta (b,\bar{a})]\, \, iff\, \, (\exists i<i^{*})(b=b_{i}\land \bar{a}=\overline{a^{i}})]\]

\end{lem}
\begin{proof}
Define \( \tau =\{c^{*}_{0},...,c^{*}_{n},c_{1},...,c_{n},s_{0},s_{1},f_{1},...,f_{n}\} \).
For each \( n\geq t\geq 0 \) define the formulas: \[
\theta _{t}(x,\bar{y}):=y_{0}=c_{0}\land ...\land y_{t}=c_{t}\land y_{t+1}=x\land y_{t+2}=f_{t+2}(x)\land ...\land y_{n}=f_{n}(x)\]
 \[
\theta (x,\bar{y}):=s_{1}(x)\land \bigwedge _{n\geq t\geq 0}[s_{0}(c^{*}_{t})\rightarrow \theta _{t}(x,\bar{y})]\]
 Now let \( U\in \frak {U} \) and assume \( \left\langle \overline{a^{i}}\in {}^{n}U:i<i^{*}\right\rangle  \)
is a delta system, this means we have some \( n\geq t^{*}\geq 0 \),
such that: \( |\{a^{i}_{t}:i<i^{*}\}|=1 \) for all \( 1\leq t\leq t^{*} \),
and \( |\{a^{i}_{t}:i<i^{*}\}|=i^{*} \) for all \( n\geq t>t^{*} \).
We can now define \( M \): \\
\( c_{0}^{*^{M}}...c_{n}^{*^{M}} \) are some distinct elements of
\( U \) (we assume \( |U|>n \)).\\
\( c_{t}^{M}=a^{1}_{t} \) (for \( 1\leq t\leq t^{*} \) and assuming
\( t^{*}>0 \) over-wise the definition of \( c^{M}_{t} \) is insignificant).\\
\( s_{0}^{M}:=\{c^{*^{M}}_{t^{*}}\} \).\\
\( s_{1}^{M}:=\{a^{i}_{t^{*}+1}:i<i^{*}\} \) (assuming \( t^{*}<n \)
over-wise define \( s^{M}_{1} \) to be some singleton).\\
\( f_{t}^{M}:=\{(a^{i}_{t^{*}+1},a^{i}_{t}):i<i^{*}\} \) (for \( t^{*}+1<t\leq n \)
and assuming \( t^{*}+1<n \) over-wise the definition of \( f^{M}_{t} \)
is insignificant).\\
Note that \( f^{M}_{t} \) are one to one functions in the relevant
cases. In conclusion we define \( \left\langle b_{i}=a^{i}_{t^{*}+1}:i<i^{*}\right\rangle  \)
(again we assume \( t^{*}<n \) over-wise we define \( \left\langle b_{i}\in U:i<i^{*}\right\rangle  \)
to be some constant sequence). So by our definitions we get \( M\models \theta _{t^{*}}(b_{i},\overline{a^{i}}) \)
for all \( i<i^{*} \). Moreover if \( M\models \theta _{t^{*}}(b,\bar{a}) \)
then there exists \( i<i^{*} \) such that \( b=b_{i} \) and \( \bar{a}=\overline{a^{i}} \).
Hence \( \theta ,M \) and \( \left\langle b_{i}:i<i^{*}\right\rangle  \)
are as needed.
\end{proof}
We now show that it is impossible to interpret large order relation
on \( \frak {U} \). 

\begin{lem}
\label{k_3^*}Let \( \tau _{0} \) be a simple vocabulary, and \( \varphi (\bar{x},\bar{y}) \)
a formula in \( \tau _{0}\cup \{r\} \) (not assuming \( lg(\bar{x})=lg(\bar{y}) \)).
Then there exists a natural number \( k^{*}=k_{3}^{*}=k_{3}^{*}(\varphi ) \)
such that for every \( \frak {M} \) a simple expansion of \( \frak {R} \)
for \( \tau _{0}\cup \{r\} \), and for all \( U\in \frak {U} \),
it is impossible to find sequences \( \left\langle \overline{a_{i}}\in {}^{lg(\bar{x})}U:i<k^{*}\right\rangle  \)
and \( \left\langle \overline{b_{j}}\in {}^{lg(\bar{y})}U:j<k^{*}\right\rangle  \)
such that:\[
(\forall i,j<k^{*})[\frak {M}[U]\models \varphi (\overline{a_{i}},\overline{b_{j}})\Longleftrightarrow i\leq j]\]

\end{lem}
\begin{proof}
Let \( \varphi (\bar{x},\bar{y}) \) and \( \tau _{0} \) be as described.
For \( i\in \{1,2\} \) let \( \tau _{i},\theta _{i} \) be the vocabulary
and formula used to code delta systems for \( n=lg(\bar{x}) \) and
\( n=lg(\bar{y}) \) respectively (i.e. those from the previous lemma).
Add to \( \tau _{0} \) new \( 1 \)-place relation symbol and function
symbol, \( s^{*},f^{*} \). In the vocabulary \( \tau =\tau _{0}\cup \tau _{1}\cup \tau _{2}\cup \{s^{*},f^{*}\}\cup \{r\} \)
(without loss of generality the unions are disjoint) define the formula:\begin{eqnarray*}
\phi (v,v'):=(\lnot \exists u)s^{*}(u)\land (\forall \bar{x},\overline{x'},\overline{y},\overline{y'})\{[\theta _{1}(u,\bar{x})\land (\theta _{1}(f^{*}(u),\overline{x'}) &  & \\
\land \theta _{2}(v,\overline{y})\land \theta _{2}(v',\overline{y'})]\rightarrow [\varphi (\overline{x},\overline{y})\equiv \varphi (\overline{x},\overline{y'})\land \varphi (\overline{x'},\overline{y})\equiv \lnot \varphi (\overline{x'},\overline{y'})]\} &  & 
\end{eqnarray*}

which will interpret a large equivalence relation. For all \( m,n\in \omega  \)
Let \( Delta(n,m) \) note the minimal number \( d \) such that every
sequence of \( d \) \( n \)-tuples has a subsequence of length \( m \)
which is a delta system. We can now define \( k^{*}_{3}(\varphi ) \):\[
k^{*}=k_{3}^{*}(\varphi ):=Delta(lg(\bar{x}),Delta(lg(\bar{y}),(k_{2}^{*}(\phi ))^{2}))\]

Seeking contradiction we assume that there exist \( \frak {M}_{0} \)
a model for \( \tau _{0} \) on \( \frak {U} \), \( U\in \frak {U} \)
and sequences as in the lemma. By the definition of \( k^{*} \) there
exist subsequences of length \( (k^{*}_{2}(\phi ))^{2} \), which
are delta systems. Note \( k_{2}:=k^{*}_{2}(\phi ) \). without loss
of generality we assume these subsequences are: \( \left\langle \overline{a_{i}}\in {}^{lg(\bar{x})}U:i<(k_{2})^{2}\right\rangle  \)
and \( \left\langle \overline{b_{j}}\in {}^{lg(\bar{y})}U:j<(k_{2})^{2}\right\rangle  \).
Let \( M_{1},M_{2} \), \( \left\langle a_{i}\in U:i<(k_{2})^{2}\right\rangle  \)
and \( \left\langle b_{j}\in U:j<(k_{2})^{2}\right\rangle  \) be
the models and sequences used to code \( \left\langle \overline{a_{i}}:i<(k_{2})^{2}\right\rangle  \)
and \( \left\langle \overline{b_{j}}:j<(k_{2})^{2}\right\rangle  \)
(see \ref{delta}). We define \( M \) a model for \( \tau  \) on
\( U \): \\
For each \( i\in \{0,1,2\} \): \( M|\tau _{i}:=M_{i} \).\\
\( s^{*^{M}}:=\{a_{j\cdot k_{2}}:j\in \{0,1,2,...,k_{2}-1\}\} \).\\
\( f^{*^{M}}:=\{(a_{j\cdot k_{2}},a_{((j+1)mod(k_{2}))\cdot k_{2}}):j\in \{0,1,2,...,k_{2}-1\}\} \).

Note that if \( \pi  \) is the permutation of \( \{\overline{a_{i\cdot k}}\in {}^{lg(\bar{x})}U:j<k_{2}\} \)
defined by \( \pi (\overline{a_{j\cdot k_{2}}})=a_{((j+1)mod(k_{2}))\cdot k_{2}} \),
then the formula:\[
\phi '(\bar{y},\overline{y'}):=(\lnot \exists \bar{x}\in \{\overline{a_{j\cdot k_{2}}}\in :j<k_{2}\})[\varphi (\overline{x},\overline{y})\equiv \varphi (\overline{x},\overline{y'})\land \varphi (\pi (\overline{x}),\overline{y})\land \lnot \varphi (\pi (\overline{x}),\overline{y'})]\]
 interprets in \( M \) a \( k_{2} \)-big equivalence relation on
\( \{\overline{b_{j}}:j<(k_{2})^{2}\} \) namely the relation \( \{(\overline{b_{i}},\overline{b_{j}}):i,j\in (k_{2})^{2},\exists l\in \{0,...,k_{2}-1\}\, s.t.\, i,j\in [l,l+1,...,l+k_{2})\} \)
. Hence by the properties of \( \theta _{1} \) and \( \theta _{2} \),
the formula \( \phi (v,v') \) interprets a \( k_{2} \)-big equivalence
relation on \( \{b_{j}:j<(k_{2})^{2}\} \) which is a contradiction.
\end{proof}
we need one more lemma before we can prove the main theorem.

\begin{lem}
\label{k_4^*}Let \( \tau  \) be a simple vocabulary and \( \varphi (x,y,\bar{z}) \)
a formula in \( \tau \cup {r} \). Then there exist a natural number
\( k^{*}=k_{4}^{*}=k_{4}^{*}(\varphi ) \) such that for every \( \frak {M} \)
a simple expansion of \( \frak {R} \) for \( \tau \cup \{r\} \)
on \( \frak {U} \) and for all \( U\in \frak {U} \), it is impossible
to find: \( \overline{c^{l}}\in {}^{lg(\bar{z})}U \) for each \( l<k^{*} \)
and sequences \( \left\langle a^{l}_{i}\in U:i<k^{*}\right\rangle  \)
and \( \left\langle b^{l}_{j}\in U:j<k^{*}\right\rangle  \) such
that:

\begin{itemize}

\item[$(\alpha)$]

For all \( l,i,j<k^{*} \), \( \frak {M}[U]\models \varphi (a^{l}_{i},b^{l}_{j},\overline{c^{l}})\, \, iff\, \, i=j \).

\item[$(\beta)$]

For all \( l_{1}<l<k^{*} \), the truth value of \( \varphi (a_{i}^{l_{1}},b_{j}^{l_{1}},\overline{c^{l}}) \)
in \( \frak {M}[U] \) is independent of \( i,j<k^{*} \).

\end{itemize}
\end{lem}
\begin{proof}
Note \( lg(\bar{z})=n \). Let \( \tau ' \) and \( \theta (x,\bar{y}) \)
be the vocabulary and formula we get by applying lemma \ref{delta}
to \( n \). Define a simple vocabulary \( \tau ^{*}:=\tau \cup \tau '\cup \{s_{1},s_{2},f\} \),
and formulas in \( \tau ^{*} \):\[
\psi _{1}(x,y):=s_{1}(x)\wedge s_{2}(y)\wedge (\forall \bar{z})[\theta (y,\bar{z})\rightarrow \varphi (x,f(x),\bar{z})]\]
 \[
\psi _{2}(x,x'):=s_{1}(x)\wedge s_{1}(x')\wedge (\forall y\forall \bar{z})(s_{2}(y)\wedge \theta (y,\bar{z}))\rightarrow [(\varphi (x,(f(x),\bar{z})\equiv \varphi (x',f(x'),\bar{z})]\]
 Put \( k^{*}:=Max\{k^{*}_{3}(\psi _{1}),k^{*}_{2}(\psi _{2})\}+1 \).
Let \( U\in \frak {U} \) and \( M \) some simple expansion of \( \frak {R}[U] \)
for \( \tau \cup \{r\} \). seeking contradiction assume that for
all \( l<k^{*} \) there exist \( \overline{c^{l}}\in {}^{lg(\bar{z})}U \)
and sequences \( \left\langle a^{l}_{i}\in U:i<k^{*}\right\rangle  \)
and \( \left\langle b^{l}_{j}\in U:j<k^{*}\right\rangle  \) satisfying
\( (\alpha ) \) and \( (\beta ) \). By increasing \( k^{*} \) to
\( Delta(n,k^{*}) \) we may assume that \( \left\langle \overline{c^{l}:}l<k^{*}\right\rangle  \)
is a delta system. Hence the truth value of the sentences \( a_{i_{1}}^{l_{1}}=a_{i_{2}}^{l_{2}} \)
and \( \varphi (a_{i}^{l_{1}},b_{j}^{l_{1}},\overline{c^{l}}) \)
depends only on the order type of the indexes (again by increasing
\( k^{*} \) and using Ramsey theorem). Moreover without loss of generality
we may assume that: \( a_{i_{1}}^{l_{1}}=a_{i_{2}}^{l_{2}}\Rightarrow (l_{1}=l_{2})\wedge (i_{1}=i_{2}) \).
This is true because we can increase \( k^{*} \) to \( (k^{*})^{2}+1 \)
and choose the sub sequence \( \left\langle a^{l}_{i}:k^{*}\cdot l\leq i<k^{*}\cdot (l+1),i\neq l\right\rangle  \)
(as \( \left\langle a^{l}_{i}:i<k^{*}\right\rangle  \)). Now if \( a_{i_{1}}^{l_{1}}=a_{i_{2}}^{l_{2}} \)
then by our assumption we have \( a_{i_{1}+1}^{l_{1}}=a_{i_{2}}^{l_{2}} \)
or \( a_{i_{1-1}}^{l_{1}}=a_{i_{2}}^{l_{2}} \), hence we get \( a_{i_{1}+1}^{l_{1}}=a_{i_{1}}^{l_{1}} \)
or \( a_{i_{1}-1}^{l_{1}}=a_{i_{1}}^{l_{1}} \) contradicting \( (\alpha ) \).
Using the same argument we may assume that \( b_{i_{1}}^{l_{1}}=b_{i_{2}}^{l_{2}}\Rightarrow (l_{1}=l_{2})\wedge (i_{1}=i_{2}) \).
Now using \( (\beta ) \) exactly one of the following conditions
hold: ever \( l_{1}<l<k^{*}\Rightarrow M\models \varphi (a_{i}^{l_{1}},b_{j}^{l_{1}},\overline{c^{l}}) \)
or \( l_{1}<l<k^{*}\Rightarrow M\models \lnot \varphi (a_{i}^{l_{1}},b_{j}^{l_{1}},\overline{c^{l}}) \).
We will deal with the first case (the second can be dealt with similarly).
We have three cases:
\begin{enumerate}
\item There exist \( i^{*}\neq j^{*}<k^{*} \) such that for all \( l<l_{1}<k^{*} \)
we have: \( M\models \lnot \varphi (a_{i^{*}}^{l_{1}},b_{j^{*}}^{l_{1}},\overline{c^{l}}) \).
\item There exists \( \pi  \) a permutation of \( \{0,...k^{*}-1\} \)
without fixed points so that for all \( l<l_{1}<k^{*} \) and for
all \( j<k^{*} \) we have: \( M\models \varphi (a_{j}^{l_{1}},b_{\pi (j)}^{l_{1}},\overline{c^{l}}) \).
\item Never (1) nor (2) hold.
\end{enumerate}
As stated above \( \tau ' \) and \( \theta (x,\bar{y}) \) are the
vocabulary and formula we get by applying lemma \ref{delta} to \( n \).
Let \( M' \) and \( \left\langle c^{l}\in U:l<k^{*}\right\rangle  \)
be the model and sequence we get by applying that lemma (in \( U \))
to \( \left\langle \overline{c^{l}:}l<k^{*}\right\rangle  \). For
each of the cases (1)-(3) we define \( M^{*} \) a simple expansion
of \( M \) for \( \tau ^{*}\cup \{r\} \) and get a contradiction.
In each case \( M^{*}|\tau ':=M' \). The interpretation of \( s_{1},s_{2} \)
and \( f \) will be given for each case separately:

Case (1): Define \( s_{1}^{M^{*}}:=\{c^{l}:l<k^{*}\} \), \( s_{2}^{M^{*}}:=\{a^{l}_{i^{*}}:l<k^{*}\} \)
and \( f^{M^{*}}:=\{(a^{l}_{i^{*}},b^{l}_{j^{*}}:l<k^{*}\} \). Then
we have \( l>l_{1}\Leftrightarrow M^{*}\models \varphi (a^{l_{1}}_{i^{*}},f(a^{l_{1}}_{i^{*}}),\overline{c^{l}}) \)
and hence the formula \( \psi _{1}(x,x') \) interprets in \( M^{*} \)
an order relation (in the sense of \ref{k_3^*}) on \( \{c^{l}:l<k^{*}\}\times \{a^{l}_{i^{*}}:l<k^{*}\} \).
This is a contradiction as \( k^{*} \) is larger than \( k^{*}_{3}(\psi _{1}) \). 

Case (2) Define \( s_{1}^{M^{*}}:=\{c^{l}:l<k^{*}\} \), \( s_{2}^{M^{*}}:=\{a_{i}^{l}:i,l<k^{*}\} \)
and \( f^{M^{*}}:=\{(a_{i}^{l},b^{l}_{\pi (i)}:i,l<k^{*}\} \). Then
we have \( l\neq l_{1}\Leftrightarrow M^{*}\models \varphi (a_{i}^{l_{1}},f(a_{i}^{l_{1}}),\overline{c^{l}}) \)
and hence the formula \( \psi _{2}(x,x') \) interprets in \( M^{*} \)
the relation \( \{(a_{i_{1}}^{l_{1}},a^{l_{2}}_{i_{2}}):l_{1}=l_{2}\wedge i_{1},i_{2}<k^{*}\} \)
which is \( k^{*} \)-big. This is a contradiction as \( k^{*} \)
is larger than \( k^{*}_{2}(\psi _{2}) \). 

Case (3): In this case we assume that (in advance) we chose \( (k^{*}+2)\cdot k^{*} \)
instead of \( k^{*} \). (So in cases (1) and (2) it is enough to
find subsets of size \( k^{*} \) with the desired properties). Look
at \( \left\langle c^{(l+2)\cdot k^{*}}:l<k^{*}\right\rangle  \),
and the sequences \( \left\langle a^{(l+2)\cdot k^{*}}_{j\cdot k^{*}}:j<k^{*}\right\rangle  \)
and \( \left\langle b^{(l+2)\cdot k^{*}}_{j\cdot k^{*}}:j<k^{*}\right\rangle  \)
for \( l<k^{*} \). Since (1) does not hold for these sequences we
get (choosing \( i^{*}=0 \) and \( j^{*}=1 \)) that there exist
\( l^{*}<l_{1}^{*} \) such that \( \varphi (a^{(l_{1}^{*}+2)\cdot k^{*}}_{0},b^{(l_{1}^{*}+2)\cdot k^{*}}_{k^{*}},c^{(l^{*}+2)\cdot k^{*}}) \).
In the same way (choosing \( i^{*}=1 \) and \( j^{*}=0 \)) we get
that there exist \( l^{**}<l_{1}^{**} \) such that \( \varphi (a^{(l_{1}^{**}+2)\cdot k^{*}}_{k^{*}},b^{(l_{1}^{**}+2)\cdot k^{*}}_{0},c^{l^{**}}) \).
Now look at \( \left\langle c^{2k^{*}+l}:l<k^{*}\right\rangle  \)
and the sequences \( \left\langle a^{2k^{*}+l}_{j}:j<k^{*}\right\rangle  \)
and \( \left\langle b^{2k^{*}+l}_{j}:j<k^{*}\right\rangle  \) for
\( l<k^{*} \). Let \( \pi  \) be a permutation of \( \{0,...,k^{*}-1\} \)
without a fixed point. We show that these sequences along with \( \pi  \),
satisfy the demands of case (2). Let \( j<k^{*} \) and \( l<l_{1}<k^{*} \).
If \( j<\pi (j) \) then \( j<\pi (j)<2k^{*}+l<2k^{*}+l_{1} \) and
\( 0<k^{*}<(l^{*}+2)\cdot k^{*}<(l_{1}^{*}+2)\cdot k^{*} \). Since
the truth value of \( \varphi  \) depends only on the order type
of the indexes we get \( \varphi (a^{2k^{*}+l_{1}}_{j},b^{2k^{*}+l_{1}}_{\pi (j)},c^{2k^{*}+l}) \)
(as \( \varphi (a^{(l_{1}^{*}+2)\cdot k^{*}}_{0},b^{(l_{1}^{*}+2)\cdot k^{*}}_{k^{*}},c^{(l^{*}+2)\cdot k^{*}}) \)).
If \( \pi (j)<j \) we get the same result, only now we use the 4-tuple
\( 0<k^{*}<(l^{**}+2)\cdot k^{*}<(l_{1}^{**}+2)\cdot k^{*} \). In
both cases we have \( \varphi (a^{2k^{*}+l_{1}}_{j},b^{2k^{*}+l_{1}}_{\pi (j)},c^{2k^{*}+l}) \)
as needed in (2). So case (3) can not hold.
\end{proof}
We are now ready to prove theorem \ref{simple} in the general case.
We prove:

\begin{thm}
\label{int R n>2}Let \( \frak {R} \) be as in \ref{Motation}.4
that satisfies \ref{assumption}. Then there exist a simple vocabulary
\( \sigma  \), \( \varphi (\bar{x}) \) a formula in \( \sigma  \)
with \( lg(\bar{x})=n(\frak {R}) \), and \( \frak {N} \) a simple
model for \( \sigma  \) on \( \frak {U} \). Such that for all \( U\in \frak {U} \)
and \( \bar{a}\in {}^{n(\frak {R})}U \):\[
\frak {N}[U]\models \varphi (\bar{a})\Longleftrightarrow (U,\frak {R}[U])\models r(\bar{a})\]
 
\end{thm}
\begin{proof}
We prove the theorem by induction on \( n(\frak {R}) \). The cases
\( n(\frak {R})=0 \) and \( n(\frak {R})=1 \) are trivial. the case
\( n(\frak {R})=2 \) was proved in \ref{n=3D2}.

Before we turn to the proof of the induction step we pay attention
to the following fact. Let \( \frak {R}' \) be as in \ref{Motation}.4.
We say that {}``\( \frak {R}' \) is definable from \( \frak {R} \)
by a simple expansion{}'' if there exist a simple vocabulary \( \tau  \),
a simple expansion \( \frak {M} \) of \( \frak {R} \) for \( \tau \cup \{r\} \)
and a formula \( \varphi (x_{0},...,x_{n(\frak {R}')-1}) \) in \( \tau \cup \{r\} \)
such that for all \( U\in \frak {U} \) and \( \bar{a}\in {}^{n(\frak {R}')}U \)
we have \( \frak {R}'[U](\bar{a}) \) iff \( \frak {M}[U]\models \varphi (\bar{a}) \).
Note that if \( \frak {R}' \) is definable from \( \frak {R} \)
by a simple expansion then \( \frak {R}' \) also satisfies assumption
\ref{assumption}. (or else \( \frak {R} \) does not satisfy the
assumption for we can define a big equivalence relation from \( \frak {R} \)
using \( \varphi  \) and the model \( \frak {M} \)). If \( \frak {R}' \)
is definable from \( \frak {R} \) by a simple expansion and \( n(\frak {R}')<n(\frak {R}) \)
then by the induction hypothesis there exist: \( \sigma _{0} \) a
simple vocabulary, \( \varphi _{0}(\bar{x}) \) a formula in \( \sigma _{0} \)
with \( g(\bar{x})=n(\frak {R}') \), and \( \frak {N}_{0} \) a simple
model for \( \sigma _{0} \) on \( \frak {U} \). Such that for all
\( U\in \frak {U} \) and \( \bar{a}\in {}^{n(\frak {R}')}U \):\[
\frak {N}_{0}[U]\models \varphi _{0}(\bar{a})\Longleftrightarrow \frak {R}'[U](\bar{a})\]
 in that case we will say that \( \frak {R}' \) satisfies the induction
hypothesis and that \( \sigma _{0},\varphi _{0} \) and \( \frak {N}_{0} \)
interprets it.

We now assume \( n(\frak {R})=n+1>2 \). We prove this case in two
stages. In the first stage we show that we can interpret the relation
\( xS_{\Delta ,\frak {M}}^{n}\bar{y} \), so we prove: 
\begin{lem}
\label{int s^n}Let \( \Delta ,\tau  \) satisfy assumption \ref{assumption2},
and let \( \frak {M} \) be a simple expansion of \( \frak {R} \)
for \( \tau \cup \{r\} \) on \( \frak {U} \). Then there exist:
\begin{itemize}
\item A simple vocabulary \( \sigma _{0} \) (\( r\not \in \sigma _{0} \)).
\item \( \varphi _{0}(x,\bar{y}) \) a formula in \( \sigma _{0} \) (\( lg(\bar{y})=n \)).
\item \( \frak {N}_{0} \) a simple model for \( \sigma _{0} \) on \( \frak {U} \).
\end{itemize}
Such that for all \( U\in \frak {U} \), \( a\in U \) and \( \overline{b}\in {}^{n}U \)
we have: \( \frak {N}_{0}[U]\models \varphi _{0}(a,\bar{b})\Longleftrightarrow aS^{n}_{\Delta ,\frak {M}[U]}\bar{b} \). 
\end{lem}
\begin{proof}
(of lemma \ref{int s^n}). Let \( \frak {M}^{*},\tau ^{*} \), \( \psi (x,x') \)
and \( \chi (x,y,\bar{z}) \) (where \( lg(\bar{z})=n-1 \)) be the
vocabulary, model and formulas interpreting \( E_{A}^{\Delta ,\frak {M}[U]} \)
and \( S^{n}_{\Delta ,\frak {M}[U]} \) that were defined in the proof
of the symmetry lemma (noted there by \( \frak {M}',\tau ' \)). We
also define a formula that will interpret an order relation in \( \frak {M}^{*} \):\[
\phi =\phi (\bar{x},\bar{y}):=[\chi (x_{0},x_{1},\bar{y})\equiv \chi (x_{2},x_{3},\bar{y})]\]
 where \( lg(\bar{x})=4 \) and \( lg(\bar{y})=n-1 \). In the vocabulary
\( \tau ^{*} \) we define a set of formulas: \( \Delta ^{*}:=\Delta \cup \{\chi ,\phi \} \).
For gravity we write \( M:=\frak {M}[U] \), \( M^{*}:=\frak {M}^{*}[U] \),
\( N_{0}:=\frak {N}_{0}[U] \) and similarly for other models, where
\( U\in \frak {U} \) is understood from the context. Next we define
some constants that we will use in the proof:
\begin{enumerate}
\item \( m_{1}:=m_{1}(\Delta ):=Max\{k^{*}_{3}(\phi ),k^{*}_{4}(\chi )\}+1 \)
for the formulas \( \chi ,\phi  \) defined above (see \ref{k_3^*}
and \ref{k_4^*}).
\item \( m_{2}:=m_{2}(\Delta ):=(m_{1})^{2}+m_{1} \).
\item For all \( U\in \frak {U} \) choose by induction on \( m_{2}\geq l \),
\( A_{l}=A^{U}_{l}\subseteq U \) such that: 

\begin{enumerate}
\item \( A_{0}=\emptyset  \).
\item \( A_{l}\subseteq A_{l+1} \) for all \( l<m_{2} \).
\item For all \( l<m_{2} \), \( r\leq n+2\cdot m_{1} \) and a type \( p\in S_{\Delta ^{*}}^{r}(A_{l},\frak {M}^{*}[U]) \):
if \( p \) is realized in \( \frak {M}^{*}[U] \) then it is realized
already in \( A_{l+1} \).
\item For all \( l<m_{2} \), \( |A_{l+1}| \) is minimal under the properties
(a)-(c).
\end{enumerate}
\item We write \( A^{*}=A^{*^{U}}=A^{U}_{m_{2}} \). 
\item Note that under these conditions there exists a bound on \( |A^{*}| \)
depending only on \( |\Delta ^{*}|,m_{1},m_{2} \) and \( n \), so
in fact the bound depends only on \( n \) and \( |\Delta | \) and
we can calculate it in the beginning of the proof. We note this bound
by \( m_{3} \). We do not calculate the value of \( m_{3} \) but
note that it increases super-exponentially as a function of \( |\Delta | \).
\item \( m_{4}:=m_{4}(\Delta ):=l^{*}(\Delta _{1}(\Delta ))+l^{*}(\Delta _{2}(\Delta ))\cdot m_{1} \).
(see \ref{bound} and \ref{symetry}).
\item \( m_{5}:=m_{5}(\Delta )=2\cdot m_{4}+m_{3}+n+2 \).
\end{enumerate}
Note by \( \frak {S}=\frak {S}^{n}_{\Delta ,\frak {M}} \) the \( n+1 \)-place
relation on \( \frak {U} \) defined by \( \frak {S}[U]:=\{(x,y,\bar{z})\in {}^{n+1}U:xS^{n}_{\Delta ,\frak {M}[U]}y\bar{z}\} \).
(We keep using the existing notation and write \( x\frak {S}[U]y\bar{z} \)
instead of \( \frak {S}[U](x,y,\bar{z}) \), or sometimes write \( xS^{n}_{\Delta ,\frak {M}[U]}y\bar{z} \)
as before). Our aim is to interpret the relation \( \frak {S} \)
by a formula in a simple model on \( \frak {U} \). First note the
following fact: Assume there exists a number \( i^{*} \) such that
for all \( U\in \frak {U} \): \( {}^{n+1}U=\bigcup _{i<i^{*}}B_{i}^{U} \).
Assume farther that for all \( i<i^{*} \) the relation \( \frak {S}_{i} \)
defined by \( \frak {S}_{i}[U]:=\frak {S}[U]\cap B_{i}^{U} \) is
interpreted by the formula \( \varphi _{i} \) and the simple model
\( \frak {N}_{i} \) for the vocabulary \( \sigma _{i} \). Then the
formula \( \bigvee _{i<i^{*}}\varphi _{i}(x,y,\bar{z}) \) in the
vocabulary \( \bigcup _{i<i^{*}}\sigma _{i} \) and the model \( \frak {N} \)
defined by \( (\forall i<i^{*})\frak {N}|\sigma _{i}=\frak {N}_{i} \)
will interpret \( \frak {S} \) as needed. We return to the proof
of the lemma. Let \( \left\langle p_{i}:i<i^{*}\right\rangle  \)
be an enumeration of all the \( \Delta ^{*} \) types of two variables
over a set of at most \( m_{3} \) parameters. Formally this means
each \( p_{i} \) is a subset of \( \Phi :=\{\varphi (x,y,u_{j_{1}},...,u_{j_{k}})\in \Delta ^{*}:k<m_{3,}j_{1},...,j_{k}\in \{0,...,m_{3}-1\}\} \).
For all \( U\in \frak {U} \) fix \( \left\langle a_{0},...,a_{l}\right\rangle  \)
some enumeration of \( A^{*^{U}} \) (of course \( l<m_{3} \)) and
we then write \( tp_{\Delta ^{*}}((a,b),A^{*},\frak {M}^{*}[U])=p_{i} \)
iff \( \frak {M}^{*}[U]\models \varphi (a,b,a_{j_{1}},...,a_{j_{k}})\Leftrightarrow \varphi (x,y,u_{j_{1}},...,u_{j_{k}})\in p_{i} \).
Note that \( i^{*} \) is uniformly bounded by \( 2^{|\Delta ^{*}|\cdot {m_{3}\choose n}} \).
For all \( i<i^{*} \) and \( U\in \frak {U} \) the \( 2 \)-place
relation on \( \frak {U} \) defined for all \( U\in \frak {U} \)
by \( \{(x,y)\in {}^{2}U:tp_{\Delta ^{*}}((x,y),A^{*^{U}},\frak {M}^{*}[U])=p_{i}\} \)
satisfies the induction hypothesis. Hence there exist a simple vocabulary
\( \sigma ^{i} \) a formula \( \varphi ^{i}(x,y) \) and \( \frak {N}^{i} \)
a simple model for \( \sigma ^{i} \) on \( \frak {U} \) such that
for all \( U\in \frak {U} \) and \( a,b\in U \):\[
\frak {N}^{i}\models \varphi ^{i}(a,b)\Longleftrightarrow tp_{\Delta ^{*}}((x,y),A^{*^{U}},M^{*}[U])=p_{i}\]
 Without loss of generality we may assume that \( \sigma ^{i} \)
has only function symbols. We use a theorem of \emph{Gaifman} about
models with a distance function (see \cite{key-4}). We get that \( \varphi ^{i}(x,y) \)
is logically equivalent to some local formula. This means for all
\( U\in \frak {U} \) the truth value of \( \varphi ^{i}(x,y) \)
in \( \frak {N}^{i}[U] \) depends only on the type of \( (x,y) \)
on the set of formulas \( \Phi ^{i}:=\bigcup _{j\in \{1,2,3\}}\Phi ^{i}_{j} \)
where:\[
\Phi _{1}^{i}:=\{f^{\varepsilon (1)}_{i}\circ f^{\varepsilon (2)}_{2}\circ ...f^{\varepsilon (t)}_{t}(x)=y:f_{1},...,f_{t}\in \sigma ^{i},\varepsilon \in {}^{t}\{1,-1\},t\leq s\}\]
 \[
\Phi ^{i}_{2}:=\{f^{\varepsilon (1)}_{i}\circ f^{\varepsilon (2)}_{2}\circ ...f^{\varepsilon (t)}_{t}(x)=x:f_{1},...,f_{t}\in \sigma ^{i},\varepsilon \in {}^{t}\{1,-1\},t\leq s\}\]
 \[
\Phi ^{i}_{3}:=\{f^{\varepsilon (1)}_{i}\circ f^{\varepsilon (2)}_{2}\circ ...f^{\varepsilon (t)}_{t}(y)=y:f_{1},...,f_{t}\in \sigma ^{i},\varepsilon \in {}^{t}\{1,-1\},t\leq s\}\]
 and \( s=s(i) \) is a natural number that depends only on \( \varphi ^{i} \).
Define for each \( j\in \{1,2,3\} \): \( \Phi _{j}:=\bigcup _{i<i^{*}}\Phi _{j}^{i} \)
and \( \Phi =\Phi _{1}\cup \Phi _{2}\cup \Phi _{3} \). Also define
\( \sigma ^{*}:=\bigcup _{i<i^{*}}\sigma ^{i} \) (w.l.o.g. the union
is disjoint) and \( \frak {N}^{*} \) is defined by \( (\forall i<i^{*})\frak {N}^{*}|\sigma ^{i}:=\frak {N}^{i} \).
Using Gaifman theorem for all \( U\in \frak {U} \) and \( a,b,a',b'\in U \)
we have \( (\bigotimes ) \):\[
tp_{\Phi }((a,b),\emptyset ,N^{*})=tp_{\Phi }((a',b'),\emptyset ,N^{*})\Rightarrow tp_{\Delta ^{*}}((a,b),A^{*},M^{*})=tp_{\Delta ^{*}}((a',b'),A^{*},M^{*})\]
 Note that \( |\Phi | \) is uniformly bounded. Moreover the bound
depends only on \( |\Delta | \) and \( n \). We treat each \( \Phi  \)
type separately this means: Let \( q \) be a type without parameters
in \( \Phi  \) (that is simply \( q\subseteq \Phi  \)). As we saw
the number of such types is bounded by \( 2^{|\Phi |} \). As we saw
in the beginning of the proof we are done if we interpret the relation
\( \frak {S}_{q} \) defined by: \( \frak {S}[U]\cap \{(x,y,\bar{z})\in {}^{n+1}U:tp_{\Phi }((x,y),\emptyset ,\frak {N}_{1}[U])=q\} \).
Clearly the relation \( \{(x,y,\bar{z})\in {}^{n+1}U:tp_{\Phi }((x,y),\emptyset ,\frak {N}_{1}[U])=q\} \)
is definable from \( \frak {N}^{*} \) by the formula \( \varphi _{q}(x,y,\bar{z}):=\bigwedge _{\phi \in q}\phi \wedge \bigwedge _{\phi \in \Phi \setminus q}\lnot \phi  \).
Now one of the following holds:
\begin{enumerate}
\item There exist \( \theta (x,y)\in \Phi _{1} \) such that \( \theta \in q \).
Then for all \( U\in \frak {U} \) and \( a,b\in U \) we have:\[
[tp_{\Phi }((a,b),\emptyset ,\frak {N}^{*}[U])=q]\Longrightarrow \frak {N}^{*}\models \theta (a,b)\]

\item For all \( \theta (x,y)\in \Phi _{1} \), \( \theta \not \in q \).
Then for all \( U\in \frak {U} \) we have:\[
\{(x,y)\in {}^{2}U:tp_{\Phi }((x,y),\emptyset ,\frak {N}^{*}[U])=q\}=\{(x,y)\in A_{q}\times B_{q}:\frak {N}^{*}[U]\models \bigwedge _{\theta (x,y)\in \Phi _{1}}\lnot \theta (x,y)\}\]
 Where we define:\[
A_{q}:=\{x\in U:tp_{\Phi _{2}}(x,\emptyset ,\frak {N}^{*}[U])=q\cap \Phi _{2}\}\]
 \[
B_{q}:=\{y\in U:tp_{\Phi _{3}}(y,\emptyset ,\frak {N}^{*}[U])=q\cap \Phi _{3}\}\]

\end{enumerate}
Assume condition (1) is satisfied. Note that for all \( U\in \frak {U} \),
\( \theta (x,y) \) defines in \( \frak {N}^{*}[U] \) a (graph of
a) 1-1 function, note this function by \( f^{U} \). The relation
defined by \( \{(x,\bar{z})\in {}^{n}U:xS^{n}_{\Delta ,\frak {M}[U]}f^{U}(x)\bar{z}\} \)
is a \( n \)-place relation definable form \( \frak {R} \) by a
simple expansion (using the formula \( (\forall t)\theta (x,t)\rightarrow \chi (x,t,\bar{z}) \)).
Hence there exist a formula \( \varphi _{1}(x,\bar{z}) \), a vocabulary
\( \sigma _{1} \) and a model \( \frak {N}_{1} \) interpreting it.
Now the formula \( \theta (x,y)\wedge \varphi _{1}(x,\bar{z})\wedge \varphi _{q}(x,y,\bar{z}) \)
and the model for \( \sigma _{1}\cup \sigma ^{*} \) which is the
union of \( \frak {N}^{*} \) and \( \frak {N}_{1} \) interprets
\( \frak {S}_{q} \) as desired.

We now assume that condition (2) is satisfied. Let \( U\in \frak {U} \)
and \( \bar{c}\in {}^{n-1}U \). We ask a question: \begin{itemize}
\item[$\diamondsuit_{q,\bar{c}}^{U}$]

Does there exist for all \( B\subseteq U \) with \( |B|\leq m_{5} \)
and \( B\supseteq A^{*} \), elements \( a,b\in U\setminus B \) such
that \( aS^{n}_{\Delta ,\frak {M}[U[}b,\bar{c} \) and \( tp_{\Phi }((a,b),\emptyset ,\frak {M}^{*}[U])=q \). 

\end{itemize}

Assume that there exist \( U\in \frak {U} \) and \( \bar{c}\in {}^{n-1}U \)
such that the answer to \( \diamondsuit ^{U}_{q,\bar{c}} \) is YES.
Choose by induction on \( j\leq m_{4} \) a pair \( (a_{j},b_{j})\in {}^{2}U \)
such that:
\begin{itemize}
\item \( a_{j}S_{\Delta ,M}^{n}b_{j}\bar{c} \).
\item \( tp_{\Phi }((a_{j},b_{j}),\emptyset ,N^{*})=q \).
\item \( a_{j},b_{j}\not \in A^{*}\cup \{a_{k}:k<j\}\cup \{b_{k}:k<j\}\cup \{c_{0},...,c_{n-2}\} \).
\end{itemize}
This is possible by the definition of \( m_{5} \) and \( \diamondsuit ^{U}_{i,\bar{c}} \).
\relax From the sequence \( \left\langle a_{0},...,a_{m_{4}}\right\rangle  \)
we omit all the elements satisfying \( a_{i}S^{n-1}_{\Delta _{1}(\Delta ),\frak {M}'[U]}\bar{c} \)
where \( \Delta _{1} \) and \( \frak {M}' \) are taken from the
symmetry lemma (see \ref{symetry}). We omitted at most \( l^{*}(\Delta _{1}) \)
elements. Now note that for all \( j_{1},j_{2} \): \( a_{j_{1}}S^{n}_{\Delta ,M}b_{j_{2}}\bar{c}\Rightarrow b_{j_{2}}S^{n}_{\Delta _{2}(\Delta ),\frak {M}'[U]}a_{j_{1}}\bar{c} \).
Hence for all \( a_{j} \) (after the change) we have \( |\{b_{j}:a_{i}S^{n}_{\Delta ,M}b_{j}\bar{c}\}|\leq l^{*}(\Delta _{2}) \).
Hence we can decrease the size of the sequences by a factor of \( l^{*}(\Delta _{2}) \)
and get \( a_{i}S^{n+1}_{\Delta ,M}b_{j}\bar{c}\Leftrightarrow i=j \).
Since the bound on \( |\Phi | \) depends only on \( n,|\Delta | \)
we may assume w.l.o.g (by increasing \( m_{1} \) and using Ramsey
theorem) that the \( \Phi  \)-type in \( N^{*} \) without parameters
of \( (a_{j_{1}},b_{j_{2}}) \) depends only on the order type of
\( (j_{1},j_{2}) \). Hence we have sequences \( \left\langle a_{0},...,a_{m_{1}}\right\rangle  \)
and \( \left\langle b_{0},...,b_{m_{1}}\right\rangle  \) such that:\begin{itemize}
\item[(*)]

For all \( j_{1},j_{2}\leq m_{1} \): \( a_{j_{1}}S^{n}_{\Delta ,M}b_{j_{2}}\bar{c}\Longleftrightarrow j_{1}=j_{2} \).

\item[(**)]

For all \( j_{1},j_{2},j_{3},j_{4}\leq m_{1} \):\\
\( tp_{\{x<y,x=y\}}((j_{1},j_{2}),\emptyset ,(\mathbb {N},<))=tp_{\{x<y,x=y\}}((j_{3},j_{4}),\emptyset ,(\mathbb {N},<))\Longrightarrow  \)\\
\( tp_{\Phi }((a_{j_{1}},b_{j_{2}}),\emptyset ,N^{*})=tp_{\Phi }((a_{j_{3}},b_{j_{4}}),\emptyset ,N^{*}) \).

\end{itemize}

Now w.l.o.g we may assume that \( m_{1}\geq |\Phi _{1}| \) (otherwise
replace \( m_{1} \) by \( max\{m_{1},|\Phi _{1}|\} \) in the definition
of \( m_{4} \)). Hence there exists \( 0<j^{*}\leq m_{1} \) such
that: \( N^{*}\models \lnot \bigwedge _{\theta (x,y)\in \Phi _{1}}\theta (a_{0},b_{j^{*}}) \)
(remember each \( \theta (x,y) \) is a function). In addition by
our definition \( tp_{\Phi }((a_{0},b_{0}),\emptyset ,N^{*})=q \),
hence by condition (2) \( a_{0}\in A_{q} \). In addition we have
\( b_{j^{*}}\in B_{q} \) as \( tp_{\Phi }((a_{j^{*}},b_{j^{*}}),\emptyset ,N^{*})=q \),
and hence \( tp_{\Phi }((a_{0},b_{j^{*}}),\emptyset ,N^{*})=q \)
(see condition (2)). In the same way we get that there exists \( 0\leq j^{**}<m_{1} \)
such that \( tp_{\Phi }((a_{j^{**}},b_{m_{2}}),\emptyset ,N^{*})=q \).
So by \( (**) \) we have \( i,j<m_{1}\Rightarrow tp_{\Phi }((a_{i},b_{j}).\emptyset ,N^{*})=q \)
and by \( (\bigotimes ) \) we get:\begin{itemize}
\item[(***)]

For all \( j_{1},j_{2},j_{3},j_{4}\leq m_{1} \): \( tp_{\Delta ^{*}}((a_{j_{1}},b_{j_{2}}),A^{*},M^{*})=tp_{\Delta ^{*}}((a_{j_{3}},b_{j_{4}}),A^{*},M^{*}) \).

\end{itemize}

We now prove:\begin{itemize}
\item[$\heartsuit$]

There exists \( m^{*}<m_{2}-m_{1}=(m_{1})^{2} \) such that if \( (a',b') \)
and \( (a'',b'') \) are pairs from \( A_{m^{*}+m_{1}} \) that satisfy
the same \( \Delta ^{*} \)-type over \( A_{m^{*}} \) in \( M^{*} \),
then \( a'S^{n}_{\Delta ,M}b'\bar{c}\equiv a''S^{n}_{\Delta ,M}b''\bar{c} \).

\end{itemize}

Assume that \( \heartsuit  \) does not hold. Then for all \( m<m_{2}-m_{1} \)
let \( (a_{m}',b_{m}') \) and \( (a_{m}'',b_{m}'') \) be pairs from
\( A_{m+m_{1}} \) realizing the same \( \Delta ^{*} \)-type over
\( A_{m^{*}} \) in \( M^{*} \), and \( \lnot (a'S^{n}_{\Delta ,M}b'\bar{c}\equiv a''S^{n}_{\Delta ,M}b''\bar{c}) \).
Choose \( \overline{c^{m}}\in {}^{n}A_{m+1} \) realizing \( tp_{\Delta ^{*}}(\overline{c^{m}},A_{m},M^{*}) \)
(this is possible, see the definition of \( A_{m+1} \)). Now look
at the formula \( \phi (\bar{x},\bar{y})\in \Delta ^{*} \). If \( l_{1}<l_{2}<m_{2}-m_{1} \)
then \( (a_{l_{2}}',b_{l_{2}}') \) and \( (a_{l_{2}}'',b_{l_{2}}'') \)
realizes the same \( \Delta ^{*} \)-type over \( A_{l_{2}} \) in
\( M^{*} \). Since \( \overline{c^{l_{1}}}\subseteq A_{l_{2}} \)
(as \( l_{1}<l_{2} \)) and since \( \chi (x,y,\bar{z})\in \Delta ^{*} \)
interprets the relation \( S^{n}_{\Delta ,M} \) in \( M^{*} \),
we get that \( a_{l_{2}}'S^{n+1}_{\Delta ,M}b_{l_{2}}'\overline{c^{l_{1}}}\equiv a_{l_{2}}''S^{n+1}_{\Delta ,M}b_{l_{2}}''\overline{c^{l_{1}}} \)
hence \( M^{*}\models \phi ((a_{l_{2}}',b_{l_{2}}',a_{l_{2}}'',b_{l_{2}}''),\overline{c^{l_{1}}}) \).
On the other hand if \( m_{1}+l_{2}\leq l_{1}<m_{2} \) then by the
choice of \( (a_{l_{2}}',b_{l_{2}}') \) and \( (a_{l_{2}}'',b_{l_{2}}'') \)
as a counter example we have \( \lnot (a_{l_{2}}'S^{n}_{\Delta ,M}b_{l_{2}}'\overline{c}\equiv a_{l_{2}}''S^{n}_{\Delta ,M}b_{l_{2}}''\overline{c}) \).
But \( a_{l_{2}}',b_{l_{2}}',a_{l_{2}}'',b_{l_{2}}''\in A_{l_{2}+m_{1}}\subseteq A_{l_{1}} \)
and \( \overline{c^{l_{1}}} \) realizes the same \( \Delta ^{*} \)-type
over \( A_{l_{1}} \) as \( \bar{c} \), so by the definition of \( \phi  \)
and \( M^{*} \) we have \( M^{*}\models \lnot \phi ((a_{l_{2}}',b_{l_{2}}',a_{l_{2}}'',b_{l_{2}}''),\overline{c^{l_{1}}}) \).
Hence if we define \( \left\langle \overline{d_{l}}=(a_{l\cdot m_{1}}',b_{l\cdot m_{1}}',a_{l\cdot m_{1}}'',b_{l\cdot m_{1}}''):l<m_{1}\right\rangle  \)
and \( \left\langle \overline{e_{l}}=\overline{c^{l\cdot m_{1}}}:l<m_{1}\right\rangle  \),
then \( \phi (\bar{x},\bar{y}) \) defines an order relation in the
sense of \ref{k_3^*} on them, in contradiction to \( m_{1}>k^{*}_{3}(\phi ) \).
This completes the proof of \( \heartsuit  \).

Now let \( m^{*} \) be the one from \( \heartsuit  \). For all \( l<m_{1} \)
we choose from \( A_{m^{*}+l+1} \) the sequence \( \overline{c^{l}}\, \, \widehat{\, }\left\langle a^{l}_{j},b^{l}_{j}:j<m_{1}\right\rangle  \)
that realizes the same \( \Delta ^{*} \)-type over \( A_{m^{*}+l} \)
in \( M^{*} \) as \( \overline{c}\, \, \widehat{\, }\left\langle a_{j},b_{j}:j<m_{1}\right\rangle  \).
We will show that these sequences satisfy the demands of lemma \ref{k_4^*}
for the formula \( \chi (x,y,\bar{z}) \). This will lead to a contradiction
as \( m_{1}>k^{*}_{4}(\chi ) \). \( (\alpha ) \) follows directly
from \( (*) \) and the equality of types. For \( (\beta ) \) let
\( l_{1}<l<m_{1} \) and \( j_{1},...,j_{4}<m_{1} \). Then \( (a^{l_{1}}_{j_{1}},b^{l_{1}}_{j_{2}}) \)
and \( (a^{l_{1}}_{j_{3}},b^{l_{1}}_{j_{4}}) \) are pairs from \( A_{m^{*}+l_{1}} \)
and from \( (***) \) and the equality of types we get that these
pairs realize the same \( \Delta ^{*} \)-types over \( A_{m^{*}+l_{1}} \)
and in particular over \( A_{m^{*}} \). So by \( \heartsuit  \)
we have \( ^{l_{1}}_{j_{1}}S^{n}_{\Delta ,M}b^{l_{1}}_{j_{2}}\bar{c}\equiv a^{l_{1}}_{j_{3}}S^{n}_{\Delta ,M}b^{l_{1}}_{j_{4}}\bar{c} \).
In conclusion as \( \bar{c} \) and \( \overline{c^{l_{2}}} \) realizes
the same \( \Delta ^{*} \)-type over \( A_{m^{*}+l_{1}} \) and by
the interpretation of the formula \( \phi  \) in \( M^{*} \) we
get \( a^{l_{1}}_{j_{1}}S^{n}_{\Delta ,M}b^{l_{1}}_{j_{2}}\overline{c^{l}}\equiv a^{l_{1}}_{j_{3}}S^{n}_{\Delta ,M}b^{l_{1}}_{j_{4}}\overline{c^{l}} \)
as in \( (\beta ) \).  

We are left with the case where for all \( U\in \frak {U} \) and
\( \bar{c}\in {}^{n-1}U \) the answer to \( \diamondsuit _{q,\bar{c}}^{U} \)
is NO. In this case: For all \( U\in \frak {U} \) and \( \bar{c}\in {}^{n-1}U \)
there exist \( B=B_{\bar{c}}\subseteq U \) with size \( \leq m_{5} \)
such that there are no \( u,v\not \in B \) satisfying \( uS^{n}_{\Delta ,\frak {M}[U]}v\bar{z} \)
and \( tp_{\Phi }((u,v),\emptyset ,\frak {N}^{*}[U])=q \). Define
\( \frak {T} \) a \( n \)-place relation on \( \frak {U} \) as
follows. For all \( U\in \frak {U} \), \( a\in U \) and \( \bar{c}\in {}^{n-1}U \):
\( \frak {T}[U](a,\bar{c}) \) iff there exists \( B\subseteq U \)
of size \( \leq m_{5} \) such that: 
\begin{itemize}
\item There are no \( u,v\not \in B \) with: \( uS^{n}_{\Delta ,\frak {M}[U]}v\bar{z} \)
and \( tp_{\Phi }((u,v),\emptyset ,\frak {N}^{*}[U])=q \).
\item \( B \) is minimal under the previous demand.
\item \( a\in B \).
\end{itemize}
It is clear that \( \frak {T} \) is definable from \( \frak {R} \)
in a simple expansion, and hence satisfies the induction hypothesis.
Let \( \sigma ^{**},\varphi ^{**}(x,\bar{z}) \) and \( \frak {N}^{**} \)
be the formula, vocabulary and model that interprets \( \frak {T} \).
We define \( m_{6}:=m_{5}\cdot Delta(m_{5},3) \) and show that for
all \( U\in \frak {U} \) and \( \bar{c}\in {}^{n-1}U \):\[
|\{x\in U:\frak {N}^{**}[U]\models \varphi ^{**}(x,\bar{c})\}|\leq m_{6}\]
 Assume towards contradiction that \( U \) and \( \bar{c} \) does
not satisfy that claim. Then by the definition of \( \varphi ^{**} \)
we have a sequence \( \left\langle B_{l}\subseteq U:l<m_{7}\right\rangle  \)
such that: 
\begin{enumerate}
\item For all \( l<m_{7} \), \( m_{5}\geq |B_{l}| \).
\item For all \( l<m_{7} \), there are no \( u,v\not \in B_{l} \) s.t.
\( uS^{n}_{\Delta ,\frak {M}[U]}v\bar{c} \) and \( tp_{\Phi }((u,v),\emptyset ,\frak {N}^{*}[U])=q \).
\item For all \( l<m_{7} \), \( B_{l} \) is minimal under (1) and (2).
\item For all \( l<m_{7} \), \( B_{l}\not \subseteq \bigcup _{m<l}B_{m} \).
\item \( \bigcup _{l<m_{7}}B_{l}=\{x:\frak {N}^{**}[U]\models \varphi ^{**}(x,\bar{c})\} \).
\end{enumerate}
To get this sequence we start with a sequence of all the sets satisfying
claims (1)-(3) in some order, and omits those that do not satisfy
claim (4). claim (5) follows straight from the definition of \( \varphi ^{**} \).
Now, by (1) and the assumption we get:\[
m_{6}<|\{x:m^{**}\models \varphi ^{*}_{i}(x,\bar{c})\}|=|\cup _{l<m_{7}}B_{l}|\leq m_{7}\cdot max\{|B_{l}|:l<m_{7}\}\leq m_{7}\cdot m_{5}\]
 so we have \( m_{7}\geq m_{6}/m_{5} \). By the definition of \( m_{6} \)
and Ramsey theorem we have \( B^{*}\subseteq U \) and \( l_{1}<l_{2}<l_{3}\leq m_{7} \)
such that \( i\neq j\Rightarrow B_{l_{i}}\cap B_{l_{j}}=B^{*} \).
We prove that \( B^{*} \) satisfies (1) and (2). Since \( B^{*}\subsetneqq B_{l_{3}} \)
this will be a contradiction to the minimallity of \( B_{l_{3}} \).
Obviously \( B^{*} \) satisfies (1), to show (2) take some \( a,b\not \in B^{*} \)
then by \( i\neq j\Rightarrow B_{l_{i}}\cap B_{l_{j}}=B^{*} \) we
have \( j\in (1,2,3\} \) such that \( a,b\not \in B_{l_{j}} \) and
since \( B_{l_{j}} \) satisfies (2) we get \( \lnot aS^{n,i}_{\Delta ,M}b\bar{c} \)
or \( tp_{\Phi }((u,v),\emptyset ,\frak {N}^{*}[U])\neq q \) as needed. 

We use Gaifman theorem again on the formula \( \varphi ^{**}(x,\bar{z}) \)
(w.l.o.g \( \sigma ^{**} \) have only function symbols). We get that
for all \( U\in \frak {U} \) the truth value of \( \varphi ^{**}(x,\bar{z}) \)
in \( \frak {M}^{**}[U] \) depends only in the type without parameters
of \( (x,\bar{z}) \) in \( \frak {M}^{**}[U] \) for the set of formulas
\( \Psi :=\bigcup _{j\in \{1,2,3\}}\Psi _{j} \) where:\[
\Psi _{1}:=\{f^{\varepsilon (1)}_{i}\circ f^{\varepsilon (2)}_{2}\circ ...f^{\varepsilon (t)}_{t}(x)=z_{i}:f_{1},...,f_{t}\in \sigma ^{**},\varepsilon \in {}^{t}\{1,-1\},t\leq s,i<n-1\}\]
 \[
\Psi _{2}:=\{f^{\varepsilon (1)}_{i}\circ f^{\varepsilon (2)}_{2}\circ ...f^{\varepsilon (t)}_{t}(x)=x:f_{1},...,f_{t}\in \sigma ^{**},\varepsilon \in {}^{t}\{1,-1\},t\leq s\}\]
 \[
\Psi _{3}:=\{f^{\varepsilon (1)}_{i}\circ f^{\varepsilon (2)}_{2}\circ ...f^{\varepsilon (t)}_{t}(z_{i})=z_{j}:f_{1},...,f_{t}\in \sigma ^{**},\varepsilon \in {}^{t}\{1,-1\},t\leq s,i,j<n-1\}\]
 and \( s \) is a natural number that depends only on \( \varphi ^{**} \).
Note that \( |\Psi | \) is uniformly bounded. Again we separate into
cases according to the \( \Psi  \)-types. Let \( q_{1},q_{2} \)
be \( \Psi  \)-types without parameters (formally \( q_{1},q_{2}\subseteq \Psi  \)).
The number of such types is bounded by \( 2^{|\Psi |} \) so as we
saw it is enough to interpret the relation \( \frak {S}_{q,q_{1},q_{2}} \)
defined by:\[
\frak {S}_{q}[U]\cap \{(x,y,\bar{z})\in {}^{n+1}U:tp_{\Psi }((x,\bar{z}),\emptyset ,\frak {N}^{**}[U])=q_{1}\wedge tp_{\Psi }((y,\bar{z}),\emptyset ,\frak {N}^{**}[U])=q_{2}\}\]
 The relation \( \{(x,y,\bar{z})\in {}^{n+1}U:tp_{\Psi }((x,\bar{z}),\emptyset ,N^{**})=q_{1}\wedge tp_{\Psi }((x,y),\emptyset ,N^{**})=q_{2}\} \)
is definable in \( \frak {N}^{**} \) by a formula noted \( \varphi _{q_{1},q_{2}}(x,y,\bar{z}) \).
Now, for \( l\in \{1,2\} \) one of the following hold: 
\begin{enumerate}
\item There exists \( \theta (x,z_{i})\in \Psi _{1} \) such that \( \theta \in q_{l} \).
Then for all \( U\in \frak {U} \), \( a\in U \) and \( \bar{c}\in {}^{n-1}U \)
we have:\[
[tp_{\Psi }((a,\bar{c}),\emptyset ,\frak {N}^{**}[U])=q_{l}]\Longrightarrow \frak {N}^{**}[U]\models \theta (a,c_{i})\]

\item For all \( \theta (x,z_{i})\in \Phi _{1} \), \( \theta \not \in q_{l} \). Then for all \( U\in \frak {U} \) we have:\[
\{(x,\bar{z})\in {}^{n}U:tp_{\Psi }((x,\bar{z}),\emptyset ,\frak {N}^{**}[U])=q_{l}\}=\{(x,\bar{z})\in A'_{q_{l}}\times B'_{q_{l}}:\frak {N}^{**}[U]\models \bigwedge _{\theta (x,z_{i})\in \Psi _{1}}\lnot \theta (x,z_{i})\}\]
 where we define:\[
A'_{q_{l}}:=\{x\in U:tp_{\Psi _{2}}(x,\emptyset ,\frak {N}^{**}[U])=q_{l}\cap \Psi _{2}\}\]
 \[
B'_{q_{l}}:=\{\bar{z}\in {}^{n-1}U:tp_{\Psi _{3}}(\bar{z},\emptyset ,\frak {N}^{**}[U])=q_{l}\cap \Phi _{3}\}\]

\end{enumerate}
Assume that there exists \( l\in \{1,2\} \) such that (1) holds.
Then, as we have seen, \( \theta (x,z_{i}) \) defines in each \( \frak {N}^{**}[U] \)
(a graph of) a one to one function. Note that function by \( f^{u} \).
the relation defined by:
 \[
( {(x,y,z_{0},...,\widehat{z_{i}},...,z_{n-2})\in {}^{n}U:xS^{n}_{\Delta ,\frak {M}[U]}y,z_{0},...,z_{i-1},f^{U}(x),z_{i+1},...,z_{n-2}} )\]
satisfies the induction hypothesis. Let \( \varphi ^{1}(x,y,z_{0},...,\hat{z}_{j},...,j_{n-2}) \)
be the formula in vocabulary \( \sigma ^{1} \) that interprets this
relation in the simple model \( \frak {N}^{1} \). In the same way
we interpret the relation defined by:\[ ( {(x,y,z_{0},...,\widehat{z_{i}},...,z_{n-2})\in {}^{n}U:xS^{n}_{\Delta ,\frak {M}[U]}y,z_{0},...,z_{i-1},f^{U}(y),z_{i+1},...,z_{n-2}} ) \]
using the formula \( \varphi ^{2}(x,y,z_{0},...,\hat{z}_{j},...,z_{n-2}) \)
in the vocabulary \( \sigma ^{2} \) and the simple model \( \frak {N}^{2} \).
Now the formula \( \varphi _{q}(x,y,\bar{z})\wedge \varphi _{q_{1},q_{2}}(x,y,\bar{z})\wedge \varphi ^{l}(x,y,z_{0},...,\hat{z}_{j},...,j_{n-2}) \)
in the vocabulary \( \sigma ^{*}\cup \sigma ^{**}\cup \sigma ^{l} \)
and the union of models: \( \frak {N}^{**},\frak {N}^{*},\frak {N}^{l} \),
interprets the relation \( \frak {S}_{q,q_{1},q_{2}} \) as needed. 

Assume then that for each \( l\in \{1,2\} \) (2) holds. Seeking a
contradiction we assume that there exists \( U\in \frak {U} \) such
that: 
\begin{itemize}
\item \( \frak {S}_{q,q_{1},q_{2}}[U]\neq \emptyset  \).
\item For all \( l\in \{1,2\} \), \( |A'_{q_{l}}|>m_{6}+|\Psi _{1}| \).
\end{itemize}
So we have \( a,b\in U \) and \( \bar{c}\in {}^{n-1}U \) such that
\( \frak {S}_{q,q_{1},q_{2}}[U](a,b,\bar{c}) \). Recall that we are
assuming \( \lnot \diamondsuit ^{U}_{q,\bar{c}} \), Hence we have
\( \frak {N}^{**}\models \varphi ^{**}(a,\bar{c})\vee \varphi ^{**}(b,\bar{c}) \),
because we can choose some minimal \( B_{\bar{c}} \) (there is one
because of \( \lnot \diamondsuit ^{U}_{q,\bar{c}} \)) and then \( a,b\not \in B_{\bar{c}} \)
is contradicting the definition of \( B_{\bar{c}} \). w.l.o.g we
assume that \( \frak {N}^{**}\models \varphi ^{**}(a,\bar{c}) \).
Now \( q_{1} \) satisfies (2) so \( \bar{c}\in B'_{q_{1}} \) as
\( tp_{\Psi }((a,\bar{c}),\emptyset ,N^{**})=q_{1} \). Note that
\( |\{a'\in A'_{q_{1}}:\frak {N}^{**}[U]\not \models \bigwedge _{\theta (x,z_{i})\in \Psi _{1}}\lnot \theta (a,c_{i})\}|\leq |\Psi _{1}| \)
(again each \( \theta (x,z_{i}) \) is a function) and hence \( A'_{q_{1}} \)
has more than \( m_{6} \) (distinct) elements \( \{a_{0},...a_{m_{6}}\} \)
satisfying \( tp_{\Psi }((a_{i},\bar{c}),\emptyset ,\frak {N}^{**}[U])=q_{1} \).
But we also have \( tp_{\Psi }((a,\bar{c}),\emptyset ,\frak {N}^{**}[U])=q_{1} \)
and \( \frak {N}^{**}\models \varphi ^{**}(a,\bar{c}) \), so for
all \( 0\leq i\leq m_{6} \), \( \frak {N}^{**}\models \varphi ^{**}(a_{i},\bar{c}) \)
which is a contradiction. 

Finally we assume that there is no \( U\in \frak {U} \) satisfying
the two demands above. We then divide \( \frak {U} \) into three
parts \( \frak {U}_{i}:i\in \{1,2,3\} \) such that: 
\begin{itemize}
\item \( \frak {S}_{q,q_{1},q_{2}}[U]=\emptyset \Longleftrightarrow U\in \frak {U}_{1} \).
\item \( |A'^{U}_{q_{1}}|<m_{6}+|\Psi _{1}|\Longleftrightarrow U\in \frak {U}_{2} \).
\item \( |A'^{U}_{q_{2}}|<m_{6}+|\Psi _{1}|\Longleftrightarrow U\in \frak {U}_{2} \).
\end{itemize}
By our assumption \( \cup _{i\in \{1,2,3\}}\frak {U}_{i}=\frak {U} \).
Now for each \( i\in \{1,2,3\} \) it is easy to interpret \( \frak {S}_{q,q_{1},q_{2}} \)
restricted to \( \frak {U}_{i} \) (using the formula \( (\exists x)x\neq x \),
or by adding a bounded number of constants to the vocabulary interpreted
as the elements of \( A'^{U}_{q_{1}} \) or \( A'^{U}_{q_{2}} \)
and using the induction hypothesis). Assume then that for each \( i\in \{1,2,3\} \)
the formula \( \varphi ^{***}_{i}(x,y,\bar{z}) \) in the vocabulary
\( \sigma ^{***}_{i} \) interprets in the model \( \frak {N}^{***}_{i} \),
the relation \( \frak {S}_{q,q_{1},q_{2}} \) restricted to \( \frak {U}_{i} \).
We now define \( \sigma ^{***}=\cup _{i\in \{1,2,3\}}\sigma ^{***}_{i}\cup \{s_{1},s_{2},s_{3}\} \)
(w.l.o.g the union is disjoint), and a model \( \frak {N}^{***} \)
for \( \sigma ^{***} \), such that for each \( i\in \{1,2,3\} \):
\( (\frak {N}^{***}|\frak {U}_{i})|\sigma ^{***}_{i}:=\frak {N}_{i}^{***} \),
and for all \( U\in \frak {U} \), \( s_{i}^{\frak {N}^{***}[U]}\neq \emptyset  \)
iff \( U\in \frak {U}_{i} \). (if \( i\neq j \) the definition of
\( (\frak {N}^{***}|\frak {U}_{i})|\sigma ^{***}_{j} \) is insignificant).
Now the formula:\[
\varphi ^{***}(x,y,\bar{z}):=\bigvee _{i\in \{1,2,3\}}(\exists us_{i}(u))\longrightarrow \varphi _{i}^{***}(x,y,\bar{z})\]
 interprets \( \frak {S}_{q,q_{1},q_{2}} \) in the model \( \frak {N}^{***} \)
as required. This completes the proof of lemma \ref{int s^n}. 
\end{proof}
In the second stage of proving theorem \ref{int R n>2} we interpret
\( \frak {R} \) itself. We prove the following:

\begin{lem}
\label{Conclusion}There exist a simple vocabulary \( \sigma  \),
and a finite set \( \Phi  \) of formulas in \( \sigma  \), and a
simple model \( \frak {N} \) for \( \sigma  \) on \( \frak {U} \).
Such that for all \( U\in \frak {U} \) and \( \bar{x},\overline{x'}\in {}^{n(\frak {R})}U \)
if \( tp_{\Phi }(\bar{x},\emptyset ,\frak {N}[U])=tp_{\Phi }(\bar{x}',\emptyset ,\frak {N}[U]) \)
then \( (U,\frak {R}[U])\models r(\bar{x})\equiv r(\bar{x}') \). 
\end{lem}
\begin{proof}
Define: \( \Delta :=\{r(x_{0},...,x_{n(\frak {R})-1})\} \), and we
note the first variable by \( x \) and the last \( n \) variables
by \( \bar{y} \) (so \( \Delta :=\{r(x,\bar{y})\} \)).

We define \( \sigma  \) and \( \Phi  \) simultanios and also we
define \( \frak {N}[U] \) for some \( U\in \frak {U} \). Let \( a,a'\in U \)
and \( \overline{b},\overline{b'}\in {}^{n}U \), and assume that
\( tp_{\Phi }(a\overline{b},\emptyset ,\frak {N}[U])=tp_{\Phi }(a'\overline{b'},\emptyset ,\frak {N}[U]) \),
where \( \sigma  \), \( \Phi  \) and \( \frak {N} \) will be defined. 

Let \( \sigma _{0},\varphi _{0}(x,\bar{y}) \) and \( \frak {M}_{0} \)
be those who interpret \( xS^{n}_{\Delta ,\frak {M}[U]}\bar{y} \),
i.e. those we get from applying the previous lemma to \( \Delta  \)
(where \( \tau =\emptyset  \) and \( \frak {M}=\frak {R} \)). For
gravity we write \( M=\frak {M}[U] \), \( \frak {N}[U]=N \), \( \frak {N}_{0}[U]=N_{0} \)
and \( R=\frak {R}[U] \). We add \( \sigma _{0} \) to \( \sigma  \),
\( \varphi _{0} \) to \( \Phi  \) and demand \( N|\sigma _{0}=N_{0} \).
Now we have:\[
aS^{n}_{\Delta ,M}\overline{b}\equiv a'S^{n}_{\Delta ,M}\overline{b'}\]
 We write \( E=E^{\Delta ,M}_{A^{\Delta ,M}_{U}} \) and define \( A^{*}=A^{\Delta ,M}_{U}\cup \{x:|x/E|\leq 2\cdot k_{1}^{*}(\Delta )\} \).
for all \( \alpha \in A^{*} \) we add to \( \sigma  \) a constant
\( c_{\alpha } \) and put \( c_{\alpha }^{N}:=\alpha  \). In addition
for each equivalence class \( x/E \) we add to \( \sigma  \) a \( 1 \)-place
relation symbol \( s_{x/E} \) and put \( s_{x/E}^{N}:=x/E \). Note
that both \( |A^{*}| \) and the number of equivalence classes is
uniformly bounded. We add to \( \Phi  \) formulas of the form \( x=c \)
and \( y_{i}=c \) for each constant \( c\in \sigma  \), and formulas
of the form \( s(x) \) and \( s(y_{i}) \) for each relation symbol
\( s\in \sigma  \). Now for each constant \( c\in \sigma  \) the
relation class on \( \frak {U} \), \( \frak {R}_{c} \) defined by
\( \frak {R}_{c}[U']:=\frak {R}[U'](c^{\frak {N}[U']},\overline{y}) \)
satisfies the induction hypothesis. That is it is a class of \( n \)-place
relation not satisfying condition (1) in theorem \ref{main}. Hence
we can add to \( \sigma  \) and \( \Phi  \) the dictionaries and
formulas we get from applying the induction hypothesis to each \( \frak {R}_{c} \),
and expand \( \frak {N} \) accordingly. Assume \( a\in A^{*} \),
then (due to the formula \( x=c_{a} \)) we have \( a=a' \). Because
of the formulas we added to \( \Phi  \) for the relation \( \frak {R}_{c_{a}} \)
we have: \[
\frak {R}_{c_{a}}[U](\bar{b})\equiv \frak {R}_{c_{a}}[U](\overline{b'})\]
 This implies \( R(c_{a}^{N},\bar{b})\equiv R(c_{a'}^{N},\overline{b'}) \).
But \textbf{}since \textbf{\( c_{a}^{N}=a=a'=c_{a'}^{N} \)} we \textbf{}get
\textbf{\( R(a,\bar{b})\equiv R(a',\overline{b'}) \)}, as claimed.
This proves the cases \( a\in A^{*} \) and \( a'\in A^{*} \).  \textbf{}

Now for each \( x/E \) (where \( x\notin A^{*} \)) and \( \bar{y}\in {}^{n}U \)
we define \( t^{x/E}_{\overline{y}}\in \{\mathbb {T},\mathbb {F}\} \)
to be the truth value the formula \( r(-,\bar{y}) \) gets for the
majority of elements in \( x/E \). This means: \( t^{x/E}_{\bar{y}}=\mathbb {T} \)
iff \( |\{x':xEx'\wedge R(x',\bar{y})\}|>k^{*}_{1}(\Delta ) \). Note
that this is true as \( x\notin A^{*} \) and so \( |x/E|>2\cdot k_{1}^{*}(\Delta ) \).
We get: \[
\neg aS_{\Delta ,M}^{n}\overline{b}\Rightarrow [R(a,\overline{b})\equiv (t^{a/E}_{\overline{b}}=\mathbb {T})]\]
 and since \( \Delta  \) has only one formula we get:\[
aS_{\Delta ,M}^{n}\overline{b}\Rightarrow [R(a,\overline{b})\equiv (t^{a/E}_{\overline{b}}=\mathbb {F})]\]
 For each \( x/E \) we have a class of relations \( \frak {R}_{x/E} \)
on \( \frak {U} \) defined by \( \frak {R}_{x/E}[U']:=\{\overline{y}\in {}^{n}U':t^{x/E}_{\overline{y}}=\mathbb {T}\} \),
which satisfies the induction hypothesis. Hence we can add to \( \sigma  \)
and \( \Phi  \) the dictionaries and formulas we get form applying
the induction hypothesis to each \( \frak {R}_{x/E} \) and expand
\( \frak {N} \) accordingly. We get for all \( x\notin A^{*} \)
:\[
\frak {R}_{x/E}[U](\bar{b})\equiv \frak {R}_{x/E}[U](\overline{b'})\]
 Since \( a/E=a'/E \) (due to the formula \( s_{a/E}(x) \)), we
have \( t^{a/E}_{\bar{b}}=t^{a/E}_{\overline{b'}} \). Assume \( \neg aS^{n}_{\Delta ',M'}\bar{b} \)
(as we saw \( aS^{n}_{\Delta ',M'}\overline{b}\equiv a'S^{n}_{\Delta ',M'}\overline{b'} \))
the we have:\[
R(a,\bar{b})\Leftrightarrow (t^{a/E}_{\overline{b}}=\mathbb {T})\Leftrightarrow (t^{a/E}_{\overline{b'}}=\mathbb {T})\Leftrightarrow (t^{a'/E}_{\overline{b'}}=\mathbb {T})\Leftrightarrow R(a',\overline{b'})\]
 as claimed. If \( aS^{n}_{\Delta ',M'}\bar{b} \) then again we get:\[
R(a,\bar{b})\Leftrightarrow (t^{a/E}_{\overline{b}}=\mathbb {F})\Leftrightarrow (t^{a/E}_{\overline{b'}}=\mathbb {F})\Leftrightarrow (t^{a'/E}_{\overline{b'}}=\mathbb {F})\Leftrightarrow R(a',\overline{b'})\]
 This completes the proof of lemma \ref{Conclusion}. 
\end{proof}
\relax From the lemma it is easy to prove that \( \frak {R} \) is interpretable
by a formula in a simple model. The proof is identical to the \( 2 \)-place
case (see the proof of \ref{n=3D2}). This completes the proof of
theorem \ref{int R n>2}.
\end{proof}

\end{document}